\documentclass[3p,authoryear,12pt]{elsarticle}

\usepackage{amsmath,amssymb,epsfig,verbatim}
\usepackage{ctable}
\usepackage{eurosym}
\usepackage{multirow}
\usepackage[toc,page]{appendix}
\usepackage{rotating}

\def\comp{\raise 1pt \hbox{$\scriptstyle\circ$}}

\def\upto{{\raise 1pt \hbox{$\scriptstyle \,\nearrow\,$}}}
\def\downto{{\raise 1pt \hbox{$\scriptstyle \,\searrow\,$}}}

\journal{}

\begin{document}
\newcommand{\bolda}{\mbox{$\mathbf{a}$}}
\newcommand{\boldb}{\mbox{$\mathbf{b}$}}
\newcommand{\boldx}{\mbox{$\mathbf{x}$}}
\newcommand{\boldy}{\mbox{$\mathbf{y}$}}
\newcommand{\boldf}{\mbox{$\mathbf{f}$}}
\newcommand{\boldz}{\mbox{$\mathbf{z}$}}
\newcommand{\boldF}{\mbox{$\mathbf{F}$}}
\newcommand{\boldG}{\mbox{$\mathbf{G}$}}
\newcommand{\boldg}{\mbox{$\mathbf{g}$}}
\newcommand{\boldh}{\mbox{$\mathbf{h}$}}
\newcommand{\boldH}{\mbox{$\mathbf{H}$}}
\newcommand{\boldzero}{\mbox{$\mathbf{0}$}}
\newcommand{\Rbb}{\mbox{$\mathbb R$}}

\newcommand{\be}{\begin{equation}}
\newcommand{\ee}{\end{equation}}

\begin{frontmatter}

\title{Optimal Management of Naturally Regenerating Uneven-aged Forests}

\author[label1]{Ankur Sinha}\ead{asinha@iima.ac.in}

\author[label3]{Janne R{\"a}m{\"o}}\ead{janne.ramo@helsinki.fi}

\author[label2]{Pekka Malo}\ead{pekka.malo@aalto.fi}

\author[label2]{Markku Kallio}\ead{markku.kallio@aalto.fi}

\author[label3]{Olli Tahvonen}\ead{olli.tahvonen@helsinki.fi}

\address[label1]{Production and Quantitative Methods, 
	Indian Institute of Management Ahmedabad\\
	Vastrapur, Ahmedabad 380015, India. \vspace{4mm}}

\address[label2]{Department of Information and Service Economy, 
Aalto University School of Business\\
PO Box 21220, 00076 Aalto, Helsinki, Finland. \vspace{4mm}}

\address[label3]{Department of Forest Sciences, University of Helsinki\\ PO BOX 27, FI-00014, Helsinki, Finland.}

\begin{abstract}{A shift from even-aged forest management to uneven-aged management practices leads to a problem rather different from the existing straightforward practice that follows a rotation cycle of artificial regeneration, thinning of inferior trees and a clearcut. A lack of realistic models and methods suggesting how to manage uneven-aged stands in a way that is economically viable and ecologically sustainable creates difficulties in adopting this new management practice. To tackle this problem, we make a two-fold contribution in this paper. The first contribution is the proposal of an algorithm that is able to handle a realistic uneven-aged stand management model that is otherwise computationally tedious and intractable. The model considered in this paper is an empirically estimated size-structured ecological model for uneven-aged spruce forests. The second contribution is on the sensitivity analysis of the forest model with respect to a number of important parameters. The analysis provides us an insight into the behavior of the uneven-aged forest model.
}
\end{abstract}

\begin{keyword}
{Natural resources, Forest management, Nonlinear programing, Evolutionary algorithms.}
\end{keyword}

\end{frontmatter}

\section{Introduction}
\color{black}Optimizing the use of forest resources has hundreds of years of history. 
The infinite horizon model specified in \cite{faustmann1849}, reintroduced by 
\cite{samuelson1976} and extended in numerous papers like \cite{kao1979notes} 
and \cite{chen1980derivation} served as a cornerstone both in research 
and in practical forestry applications.~\color{black}
In its generic form the model determines optimal forest rotation, i.e., the length of optimal interval between clearcuts. However, it is seldom noticed that since this model is most suitable for plantations 
\citep{yoshimoto1998searching}, it has directed research to forests that actually cover only $7\%$ of the total world forest land area. Our study presents major progress in the line of research that serves developing the management of more natural forest \color{black}stands~\color{black} that have great potential in solving several pressing problems related to forest environment.

	The alternative to plantations is to rely on native tree species, natural regeneration and continuous forest cover, i.e., to manage forests as heterogeneous uneven-aged systems. The rationale of this model depends on tree species, but for shade-tolerant trees the economic outcome may become fully competitive because of \color{black}natural~\color{black} regeneration and more accurate targeting of cuttings to those trees that are financially mature. Additionally, managing forest resources in more natural and heterogamous state has high potential in coping with problems such as climate change \citep{field2014}, loss of biodiversity and landscape esthetics \citep{thompson2009}. Multi-criteria decision making approaches have also been used in the past to meet multiple objectives in forest management problems \citep{steuer1978interactive,nhantumbo2001goal}. Interest in continuous cover forestry is increasing in Nordic countries and UK, for e.g., in Finland it has been released from a 70 year of legislation ban from the beginning of 2014. According to surveys a major problem among forestry professionals is the suspense of the alternative system's economic viability \citep{valkonen2014}.

	While the Faustmann approach describes a chain of exactly similar even-aged cohorts, the model for more natural forests includes the internal structure of heterogeneous trees. As shown in the seminal paper by \cite{adams1974} this expands model dimensions and the development of the research has been a struggle against limitations in computing capacity. This has led researchers to develop various simplifications with the cost of losing economically and mathematically sound theoretical structure as already surveyed by \cite{getz1989}. Most studies still circumvent the problems by studying the fundamentally dynamic problem in a static setup with limited scientific progress and low practical credibility. One problem is in solving a multiple state variable infinite horizon model from any initial \color{black}stand~\color{black} state. A straightforward solution for this problem was given already in \cite{haight1990}: lengthen the planning horizon until the approach path towards a stationary state (or cycle) becomes invariant from further lengthening and consider it as an approximation of the infinite horizon solution. Given a single tree species cases this leads to solvable problems with e.g., 200 periods and 24 optimized variables per period; albeit non-linearities and non-convexities require special attention. Besides dimensionality the other problem is that the optimal solution becomes cutting \color{black}stand~\color{black} every period which does not make sense in the presence of fixed harvesting cost and the fact that too small yield is commercially invaluable. Fixed harvesting cost is taken into account in even-aged models, like \cite{tahvonen2013}, with the implication that optimal number of intermediate cuttings is between zero and five periods (25 years) depending on factors like site fertility and interest rate. However, in the even-aged problem all rotations are similar implying that the time horizon in computation is relatively short (40-150 years) and the number of combinational variables are usually six or lower. In uneven-aged models \cite{haight1990} take this into account by allowing cuttings every 20 years only. \cite{wikstrom2000} includes fixed harvesting cost and computes solutions using tabu search under the simplification that regeneration is fixed to 50 trees per 5 years period and that stand volume is not allowed to decrease below a level determined by Swedish forest legislation. He does not interpret his results on harvesting interval but it seems to vary between 5 and 20 years without any systematic pattern. In \cite{tahvonen2011} the model includes fixed harvesting cost which leads to optimal harvesting period of 15-20 years under the constraint that the interval is constant over time.

	Given these studies, the proper solution method and most general solutions for the uneven-aged management problem are still open. This is pressing in the practically most important cases where the initial forest state is a consequence of even-aged management and the problem is to solve optimal path or transition to uneven-aged management. This question has been studied in numerous works with specifications without full generality. In this paper we make a two-fold contribution. As the first contribution we develop a computational method for solving uneven-aged \color{black}stand~\color{black} management problems that is a large scale mixed integer non-linear program; and as the second contribution we provide an analysis for the uneven-aged \color{black}stand~\color{black} management model.  The solution method for handling the problem is based on the following:
\begin{enumerate}
\item A two-level approach with genetic algorithm at the upper level and continuous non-linear programing at the lower level: The approach is faster by more than an order of magnitude in terms of computation time as compared to branch-and-bound method. This supports handling of large scale uneven-aged management problems.
\item Modeling the infinite time horizon uneven-aged \color{black}stand~\color{black} management problem into a tractable problem by assuming transition and steady states: The assumption causes no loss of generality as the transition and steady state lengths are assumed to be endogenous subject to optimization. For an earlier study on forest management practices where the time horizon is divided into transition and steady states the readers may refer to \cite{salo2003}.
\end{enumerate}

A faster algorithm allowed us to perform a number of computational studies by varying the parameters in the uneven-aged \color{black}stand~\color{black} management problem. This provided us an insight into the behavior of the uneven-aged model. These insights may play a significant role in directing future research on uneven-aged management. 

The later part of the paper is structured as follows. In Section \ref{sec:model} we discuss the size-structured \color{black}stand~\color{black} model and introduce the net present value maximization problem. This is followed by the description of the algorithm in Section \ref{sec:algorithm} that is used for solving the optimization problem. Thereafter, in Section \ref{sec:results} we present the results and provide comparisons against the standard approaches that are used to solve uneven-aged stand management problems. Finally, the conclusions are provided in Section \ref{sec:conclusions}, where we also highlight the future research directions on uneven-aged \color{black}stand~\color{black} management. The paper includes appendices that provides additional computational results.

\section{\bf Size-structured Forestry model}\label{sec:model}
The model being considered in this paper is a discrete infinite time horizon model that involves two kinds of variables that are listed below:
\begin{enumerate}
\item Binary variables representing the harvesting stages, i.e., whether to harvest or not to harvest at a particular time stage.
\item Continuous variables that define the state of the forest and the extent of harvests at each time stage, among other variables.
\end{enumerate}
\begin{figure}[!ht]
\begin{center}
\epsfig{file=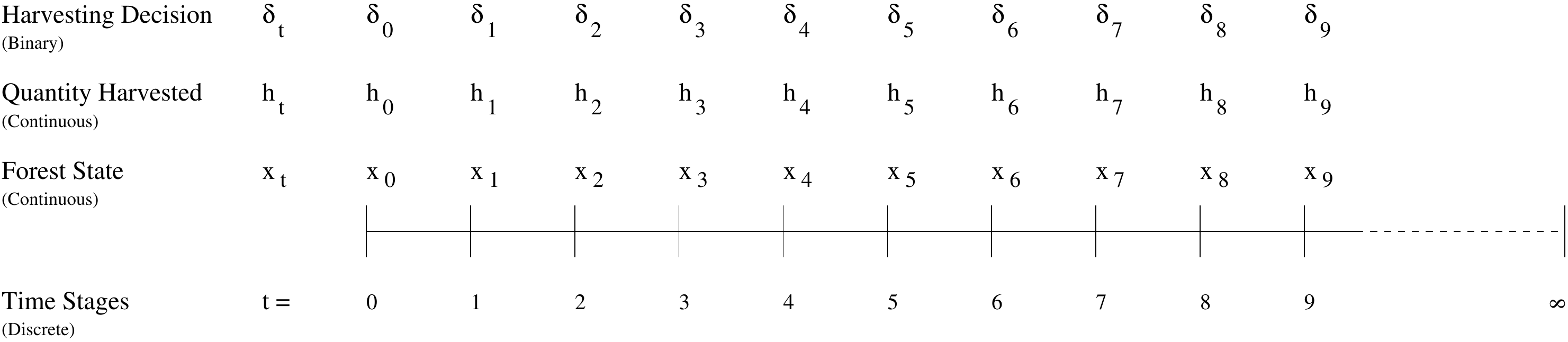,width=0.97\linewidth}
\end{center}
\caption{Forest management strategy on an infinite time horizon.}
\label{fig:forestManagementStrategy}
\end{figure}
A forest management strategy is shown in Figure~\ref{fig:forestManagementStrategy} that we want to optimize for maximum net present value ($NPV$) over an infinite horizon
in a discrete time framework. The time stages are represented as $t=0,1,2,\ldots$ on an infinite time horizon. Harvesting stages are represented by $\delta_t$ that takes values $0$ or $1$ with $1$ denoting that harvesting is done and $0$ denoting that no harvesting is done at a given time stage. The forest states and the extent of harvests are represented with vectors $x_t$ and $h_t$ respectively. We discuss forest land of one hectare. Larger forest areas require minor and  straightforward modifications which we omit. The size-structured forestry model defined in this section utilizes a number of symbols that are described in the discussions. For ease of reference we have also provided these symbols in Table \ref{tab:notations}.


\begin{table}[]
\centering
\caption{Notations used in the size-structured forestry model}
\label{tab:notations}
\begin{tabular}{@{}ll@{}}
\toprule
Important symbols &   \\ \midrule
$s$: Size class  & $t$: Time stage  \\
$x_{st}$: Trees in size class $s$ at time stage $t$              &  $x_{t}$ = ($x_{st}$) $\forall\; s$\\
$h_{st}$: Harvest in size class $s$ at time stage $t$              &  $h_{t}$ = ($x_{st}$) $\forall\; s$\\
$\delta_t$: Binary (0/1) harvesting decision at $t$
              & $\Delta$: Single time step in years\\ 
$b_s$: Basal area of tree in size class $s$ 
              & $B_t$: Total basal area per hectare \\ 
$B_{st}$: Total basal area in size class more than $s$
              & $\phi_t$: Ingrowth into size class 1 in $\Delta$ step\\
$\mu_{st}$: Mortality share of trees in size class $s$
              & $\alpha_{st}$: Share of trees growing from $s$ to $s+1$\\
$R_t$: Gross revenue & $C_t$: Cutting and hauling cost\\
$C_c$: Logging cost function& $C_h$: Hauling cost function\\
$C_f$: Fixed cost & $\beta$: Discounting factor\\ 
$v_{1s}$: Tree volume of small diameter logs in $s$ & $v_{2s}$: Tree volume of saw timber in $s$\\
$v_{s} = v_{1s}+v_{2s}$& $v = (v_s) \; \forall \; s$ \\ 
\midrule
Other symbols &   \\ \midrule
Ingrowth related: $S_1$, $S_2$, $\gamma$, $\nu$, $B^0$ & Mortality related: m   \\
Growth related: $A_1$, $A_2$, $S$, $L$ & Discounting related: $r$\\
Revenues related:  $P_1$, $P_2$, $C_1$, $C_2$&   \\
\bottomrule
\end{tabular}
\end{table}

Trees in the forest are subdivided into a finite number of size classes $s$ for $s=1,2,\dots,n$ in increasing order. Let $x_{st}$ be 
the number of trees in size class $s$ at stage $t$ and define vector $x_t=(x_{st})$. For $t=0$, $x_0$ is the given initial state of the forest.
Let vector $h_t=(h_{st})$ denote the level of harvesting at stage $t$. 
Component $h_{st}$ is the number of trees harvested in size class $s$ at time stage $t$.
For all $t$, let $\delta_t$ be a binary variable indicating whether harvesting 
takes place at stage $t$ ($\delta_t=1$) or not ($\delta_t=0$). Then a logical requirement for harvesting levels is
\be
h_t=\delta_t h_t.
\label{dh}
\ee
Before stating the optimization problem we introduce a number of endogenous auxiliary variables concerning forest dynamics and cash flow. \cite{bollandsas2008} use Norwegian National Forest Inventory data and estimate uneven-aged models for most common Nordic tree species. In our study we use their data for Norway spruce.

Given basal area $b_s$ of a tree in size class $s$, the total basal area (per hectare) 
at stage $t$ is
\be
B_t = \sum_s b_s x_{st}
\label{bt}
\ee
and the total basal area (per hectare) of trees in size classes larger than $s$ is
\be
B_{st} = \sum_{i>s} b_i x_{it}.
\label{byi}
\ee
Ingrowth $\phi_t$ of trees in step $\Delta$ into the smallest size class 1
as a function of basal area $B_t$ is
\be
\phi_t = \frac{S_1 (B_t+B^0)^{-\nu}}{1+S_2\exp(\gamma B_t)}
\label{st}
\ee
where $S_1$, $S_2$ , $B^0$, $\gamma$ and $\nu$ are positive parameters.
Mortality $\mu_{st}$ is the share of trees dying in size class $s$ in one step $\Delta$.
As a function of basal area $B_t$, it is given by 
\be
\mu_{st} = \frac{1}{1+M_s\exp(-m B_t)}
\label{mti}
\ee
where $m$ is a positive parameter, and $M_s=exp(2.492+0.02 d_s - 3.2 \times 10^{-5} d_{s}^{2})$ depends on tree diameter $d_s$ in each size class $s$.
For size class $s$, given basal area  $B_t$ and 
the basal area $B_{st}$ of trees in size classes larger than $s$, 
the share of trees growing in one step $\Delta$ from size class $s$ to $s+1$ is
\be
\alpha_{st} = \left\{\begin{array}{ll}
 G_s(S,L)-A_1 B_{st}-A_2 B_t & {\rm ~~~for~} s<n \\
 0   & {\rm ~~~for~} s=n
\end{array} \right.
\label{gti}
\ee
where $G_s$ depends on site index $S$ and latitude $L$ as follows
$$
G_s(S,L) = 0.02( 17.839+0.0476 d_s - 11.585 \times 10^{-5} d_{s}^2 + 0.906 S - 0.268 L ),
$$
while $A_1$ and  $A_2$ are positive parameters. In this notation, the share of trees remaining in size class $s$ is $1-\mu_{st}-\alpha_{st}$.

The gross revenue $R_t$ at any time step $t$ is given as
\be
R_t =  \sum_{s} h_{st} (v_{1s} p_1 + v_{2s} p_2)
\label{vt}
\ee
where $v_{1s}$ is the tree volume of small diameter logs in size class $s$, $p_1$ is the unit price of small diameter logs, $v_{2s}$ is the volume for saw timber in size class $s$ and $p_2$ is the unit price for saw timber. Cutting and hauling cost $C_t$ depend on the overall volumes in each size class $v_s=v_{1s}+v_{2s}$. Let $v=(v_s)$ be the vector of volumes, then $C_t$ is given as follows:
\be
C_t =  C_{c}(h_t,v)+C_{h}(h_t,v) + \delta_t C_{f}
\label{Ct}
\ee
where $C_c$ is the cost function for logging cost and $C_h$ is the cost function for hauling cost, and fixed cost is given by $C_f$. If $\beta$ denotes the annual discount factor per year, then the problem is 
to find $x_t\geq 0$, $h_t\geq 0$, binary variables $\delta_t$
as well as auxiliary variables $B_t$, $B_{st}$, $\phi_t$, $\mu_{st}$, $\alpha_{st}$, 
$R_t$ and $C_t$, for all $t$ and $s$, to
\be
\max\sum_{t=0}^{\infty} (R_t - C_t)\beta^{t\Delta}
\label{obj}
\ee
subject to (\ref{dh})--(\ref{Ct}) and forest dynamics state equations
\be
x_{1,t+1} = \phi_t + (1 - \mu_{1,t} - \alpha_{1,t})x_{1,t} - h_{1,t}~~~~~~~{\rm for~all}~t
\label{xt1}
\ee
\be
x_{s+1,t+1} = \alpha_{s,t} x_{s,t} + (1 - \mu_{s+1,t} -  \alpha_{s+1,t})x_{s+1,t} - h_{s+1,t}~~{\rm for~all}~t~{\rm and~}s<n
\label{xts}
\ee
\be
x_{0} = x^0
\label{x0}
\ee
The functions $C_{c}(h_t,v)$ and $C_{h}(h_t,v)$ used in the model have been defined below. The harvesting cost is specified following the estimation results in \cite{nurminen2006} assuming that in uneven-aged management cutting costs are $15\%$ higher than in even-aged clearcut operations \citep{surakka2007poimintahakkuiden}, while hauling cost is determined as in even-aged thinning operations \citep{tahvonen2011}.
\begin{align}
C_{c}(h_t,v) &= C_1 \sum_s h_{st} \big( 0.412+0.758 v_s+0.180 v_{s}^{2} \big)\label{Cc}\\
C_{h}(h_t,v) &= C_2 \big( 14.83 \delta_t + 2.272 \sum_s h_{st} v_s + 0.5348\big(\sum_s h_{st} v_s\big)^{0.7}\big)
\label{Ch}
\end{align}
where $C_1$ and $C_2$ represent the cutting cost per minute and hauling cost per minute. It is noteworthy that the hauling cost has an inherent fixed term. For ease of discussions in a later section we denote hauling cost as $C_{h}(h_t,v) = 14.83 C_2 \delta_t + C_{h}(h_t,v)'$, where the first term is the fixed term and the second term is the variable hauling cost.

The parameter values employed in our numerical illustrations are given in 
Tables~\ref{tab:data_i}-\ref{tab:data}. Table~\ref{tab:data_i} provides the size class 
dependent parameters in the model, i.e basal area ($b_s$), diameter ($d_s$), and tree 
volumes \citep{heinonen1994} of small diameter logs ($v_{1s}$) and saw timber ($v_{1s}$). 
It also provides three different initial states of the forest that have been studied in this paper. 
The three initial states $x_{0} = x^{1}, x^{2}, $ and $x^{3}$ represent an old even-aged stand, 
uneven-aged stand and young even-aged stand respectively. Table~\ref{tab:data} provides other 
parameters of the model that are independent of the size class.
\color{black}All prices and costs are given at the level of the year 2011.
Estimated costs correspond to average costs by hectare of large enough stands. 
Note also that economies of scale may be taken into account 
by varying the fixed cost.\color{black}

 Given that dimension of cash flow $R_t-C_t$ is \euro, 
basal areas $B_t$ and $B_{st}$ are m$^2$, volume $v_s$ is m$^3$,
and step size $\Delta$ is years, parameter dimensions in
Tables~\ref{tab:data_i}-\ref{tab:data} are implied by (\ref{dh})--(\ref{obj}).

\begin{table}[ht]
\caption{Tree data by size class $s$ for $n=12$ classes.}
\label{tab:data_i}
\begin{center}
\begin{tabular}{|c|r|r|r|r|r|r|r|}
\hline
\multicolumn{5}{|c|}{} & 
\multicolumn{3}{c|}{Initial States} \\
\hline
\multicolumn{1}{|c|}{$s$} &
\multicolumn{1}{c|}{$b_s$ (m$^2$)} &
\multicolumn{1}{c|}{$d_s$ (mm)} &
\multicolumn{1}{c|}{$v_{1s}$ (m$^3$)} &
\multicolumn{1}{c|}{$v_{2s}$ (m$^3$)} &
\multicolumn{1}{c|}{$x_{0}=x^{1}$} &
\multicolumn{1}{c|}{$x_{0}=x^{2}$} &
\multicolumn{1}{c|}{$x_{0}=x^{3}$} \\
\hline
1	&	0.0440	&	75	&	0.014	&	0		&	1750	&	50	&	190	\\
2	&	0.0123	&	125	&	0.067	&	0		&	0	&	25	&	162	\\
3	&	0.0241	&	175	&	0.167	&	0		&	0	&	10	&	140	\\
4	&	0.0398	&	225	&	0.081	&	0.234	&	0	&	0	&	124	\\
5	&	0.0594	&	275	&	0.065	&	0.446	&	0	&	25	&	75	\\
6	&	0.0830	&	325	&	0.060	&	0.684	&	0	&	250	&	18	\\
7	&	0.1104	&	375	&	0.050	&	0.963	&	0	&	25	&	0	\\
8	&	0.1419	&	425	&	0.050	&	1.253	&	0	&	0	&	0	\\
9	&	0.1772	&	475	&	0.043	&	1.574	&	0	&	0	&	0	\\
10	&	0.2165	&	525	&	0.039	&	1.900	&	0	&	0	&	0	\\
11	&	0.2597	&	575	&	0.033	&	2.214	&	0	&	0	&	0	\\
12	&	0.3068	&	625	&	0.031	&	2.565	&	0	&	0	&	0	\\
\hline
\end{tabular}
\end{center}
\end{table}
\begin{table}[ht]
\footnotesize
\caption{Data parameters independent of trees' size class}
\label{tab:data}
\begin{center}
\begin{tabular}{|lc|cclcclcclccl|}
\multicolumn{14}{c}{} \\
\hline 
ingrowth   	&(\ref{st})  & $S_1$   &=&147.8    & $S_2$  &=&0.5494    & $\gamma$&=&0.0180 & $\nu$   &=&0.157 \\
           	&            & $B^0$   &=&0.741    &        & &         &         & &       & &&\\
mortality  	&(\ref{mti}) & $m$   &=&0.0310   &        & &         &         & &       & &&\\
growth     	&(\ref{gti}) & $A_1$&=&0.006824 & $A_2$&=&0.000480 &         S&=&15       & L&=&60 (deg.)\\
discounting	&(\ref{obj}) & $\beta$&=&$1/(1+r)$& $r$    &=&0.03     & $\Delta$&=&5      & &&\\ 
Revenues		&(\ref{vt})  & $p_1$&=&\euro 34.07 & $p_2$    &=&\euro 58.44      & $C_1$&=&\euro 2.1      & $C_2$&=&\euro 1 \\ 
			&            & $C_f$   &=&300    &        & &         &         & &       & &&\\
\hline
\end{tabular}
\end{center}
\end{table}
All the parameters and functions have been kept fixed as suggested above. However, we have performed certain sensitivity studies by varying the parameters $r$, $C_f$ and $S$.

\section{Proposed Algorithm}\label{sec:algorithm}
Finding an optimal forest management strategy requires the optimization of binary 
as well as continuous variables, which makes the problem in the previous section 
a mixed integer non-linear programing problem. 
Mixed integer programing commonly arises 
in \color{black}optimization of forest harvesting operations~\color{black}
and researchers have used both heuristics \citep{weintraub1994heuristic} and exact 
methods \citep{goycoolea2005harvest,constantino2008new,carvajal2013imposing} 
to handle \color{black} the problems~\color{black} in different contexts. Both the approaches 
have its own advantages and disadvantages. An exact method, in such cases, would 
guarantee an optimal solution but might be computationally intractable for large 
scale problems. On the other hand a heuristic might be computationally tractable, 
but does not guarantee optimality. 

In this section, we describe a two level approach that is customized to solve uneven-aged forest management problem. The harvesting decisions are determined ($\delta_t$) at {\em level 1} of our solution procedure that fixes the values of the binary variables at each time stage to 0 or 1. Once the binary variables are fixed, at {\em level 2} we determine the corresponding optimal state of the \color{black}stand~\color{black} ($x_s = (x_{st})$) and the extent of harvests ($h_s = (h_{st})$) at each time stage. The harvesting strategies at {\em level 1} are generated using an evolutionary algorithm and their corresponding optimal continuous variables are determined at {\em level 2} by solving a non-convex optimization problem. The proposed method performs this process iteratively with an intelligent update of the harvesting strategies and leads to a near optimal solution. This section is divided into two parts: the first part provides a detailed description of a customized evolutionary procedure to optimize harvesting strategy and also describes how we convert the infinite time horizon problem into a finite time horizon problem; the second part involves a discussion about the non-convex optimization problem being solved for each harvesting strategy.

\subsection{Level 1: Evolutionary Optimization}\label{sec:level1}

In this section we describe the evolutionary optimization algorithm that uses principles from biological evolution to move towards the optimum by generating improved harvesting strategies. It is a population based approach where each member represents a harvesting strategy. The quality of a harvesting strategy is measured using {\em level 2} that returns the maximum net present value for the given strategy. We refer the net present value corresponding to each harvesting strategy as the fitness of the strategy. The technique emphasizes better harvesting strategies in the population leading to a rise in average fitness over generations (iterations). We use a genetic representation to code harvesting strategies in our algorithm. The algorithm begins with random initialization of a population representing different harvesting strategies. The strategies are then improved over generations by repetitively applying Selection, Crossover, Mutation and Replacement that is described later in this section. For few earlier studies where evolutionary computation techniques have been used in the context of forest management, the readers may refer to \color{black} \cite{bayat2013productivity}, for instance.~\color{black}

\subsubsection{Genetic Representation}
\begin{figure}[!ht]
\begin{center}
\epsfig{file=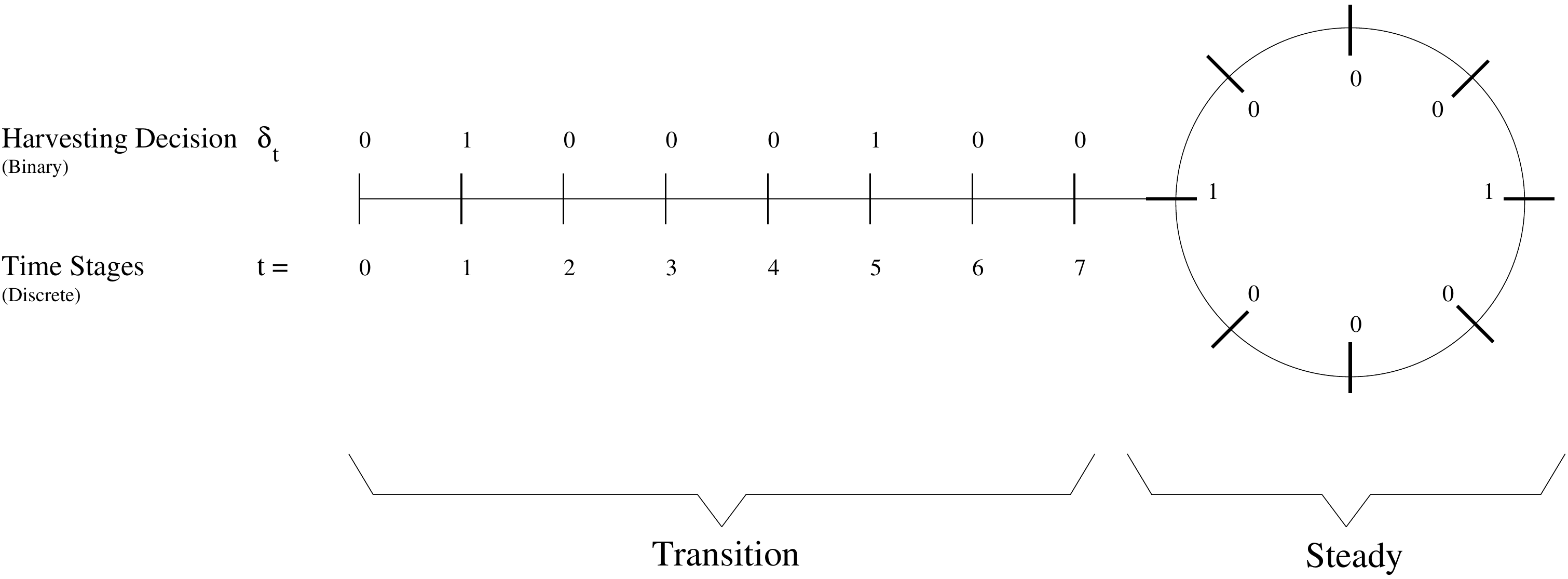,width=0.97\linewidth}
\end{center}
\caption{Forest management strategy with transition and steady states.}
\label{fig:infiniteToFinite}
\end{figure}

It is commonly observed in discrete-time dynamic systems that the variables go through a transition phase and eventually stabilize into a steady state where a fixed pattern gets repeated. Taking insights from such behavior of discrete-time dynamic systems and our prior experience in solving forest management problems \citep{salo2003}, we start with an assumption that the optimal solution to the forest management problem on an infinite horizon consists of transition period and steady state period. During the transition period the state and the control variables keep changing without following any regular pattern. However, during the steady state period both the state and control variables change in a cyclic manner. Therefore, we model an infinite time horizon with a finite set of variables as shown in Figure~\ref{fig:infiniteToFinite}. In the figure, the variables vary over time during the transition phase and once the steady state begins, the harvesting decisions, quantity harvested, as well as the states of the forest follow a repeated cycle. Such a construct allows us to explore a limited time horizon leading to a significant reduction in the search region. Using this representation, any harvesting strategy with transition and steady periods can be coded with finite binary variables. Within the algorithm the solution shown in Figure~\ref{fig:infiniteToFinite} will have a genetic representation as $\{(01000100),(10001000)\}$. The number of bits in the transition period denotes the transition period length, and the number of bits in the steady state period denotes the cycle length. The length of the transition period and steady state cycle is automatically adapted by our algorithm.

\subsubsection{Population Initialization}\label{sec:initialization}
The initial population consisting of $N$ members is generated randomly by the algorithm. The minimum and the maximum transition length ($t_{\mbox{min}},t_{\mbox{max}}$) and steady state cycle length ($s_{\mbox{min}},s_{\mbox{max}}$) are provided as input by the user. For a given member, the algorithm decides the transition length ($t_{\mbox{len}}$) by generating a random integer between $t_{\mbox{min}}$ and $t_{\mbox{max}}$. Thereafter, it generates a genetic string of length $t_{\mbox{len}}$ with 0 or 1 appearing at each location with equal probabilities. A genetic string of length $s_{\mbox{len}}$ for steady state cycle is generated in a similar manner, and then the two strings are combined to obtain a harvesting strategy for the member. This operation is repeated for every member in the population leading to a wide variety of harvesting strategies.

\subsubsection{Fitness Assignment}\label{sec:fitnessAssignment}
Each harvesting strategy or population member generated in the algorithm has to be evaluated in terms of the maximum net present value that can be obtained with the strategy over an infinite period of time. Identifying this maximal value itself is a non-convex optimization problem that we solve at {\em level 2}. The non-convex optimization problem is solved with respect to variables $x_t$ and $h_t$ leading to their optimal values corresponding to the given harvesting strategy ($\delta_t$). It is noteworthy that we need to solve a non-convex optimization task for every new harvesting strategy that we generate making the overall task computationally demanding. The non-convex optimization is discussed in detail in Subsection~\ref{sec:nonConvex}.

\subsubsection{Genetic Operators}\label{sec:geneticOperators}
\begin{figure}[t]
\begin{center}
\epsfig{file=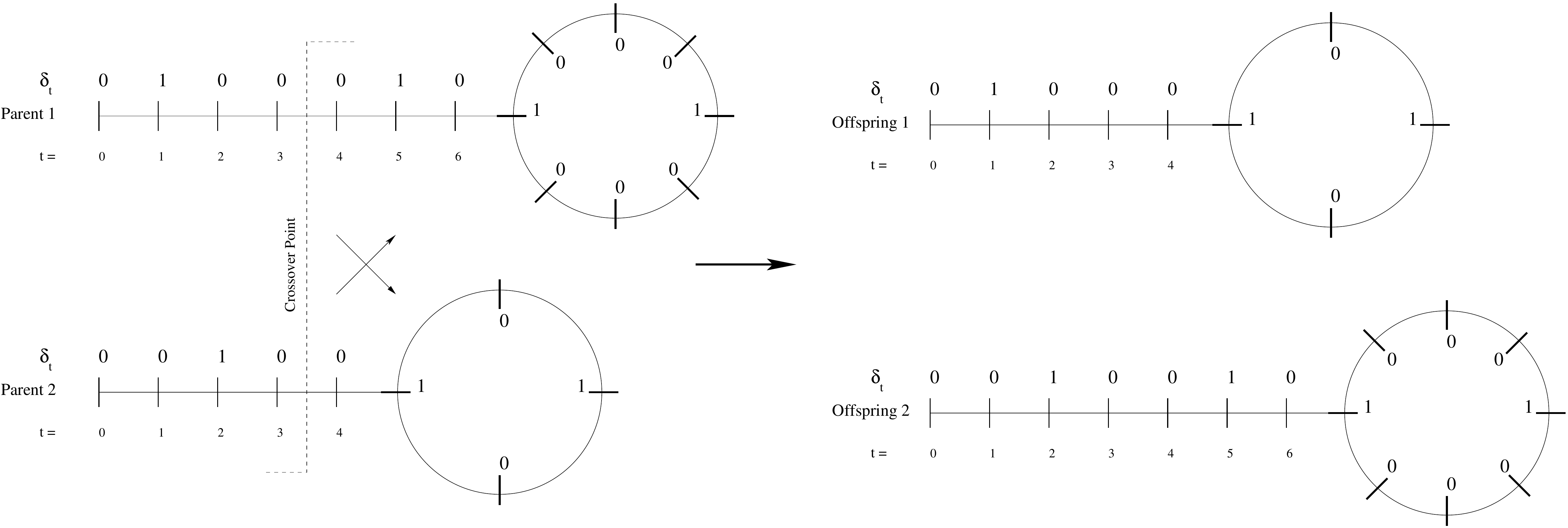,width=0.99\linewidth}
\end{center}
\caption{Crossover between two parent harvesting strategies leading to offspring strategies.}
\label{fig:crossover}
\end{figure}

The genetic operators are used in an evolutionary algorithm to generate new members using the existing population members. The existing members that are used to generate new members are referred to as parents while the newly generated members are referred to as offspring. The genetic operators consist of the steps of crossover and mutation. Crossover involves information sharing between two or more parents while mutation is performed on a single member and is helpful in maintaining diversity. 
The crossover mechanism used in this study involves two parents that lead to two offspring. The crossover operator is applied on the parents with a probability of $p_c$. In case no crossover is performed between the parents then the offspring are considered to be identical to the parents. Figure~\ref{fig:crossover} shows the crossover operation between two parents with unequal transition and steady state cycle lengths. To perform a crossover a random crossover point is chosen on the member with a smaller transition length. We fix this crossover point on both the members and swap the genetic string that follows the crossover point to generate two offspring. Next, the offspring generated using this operation undergoes mutation.

\begin{figure}[b]
\begin{center}
\epsfig{file=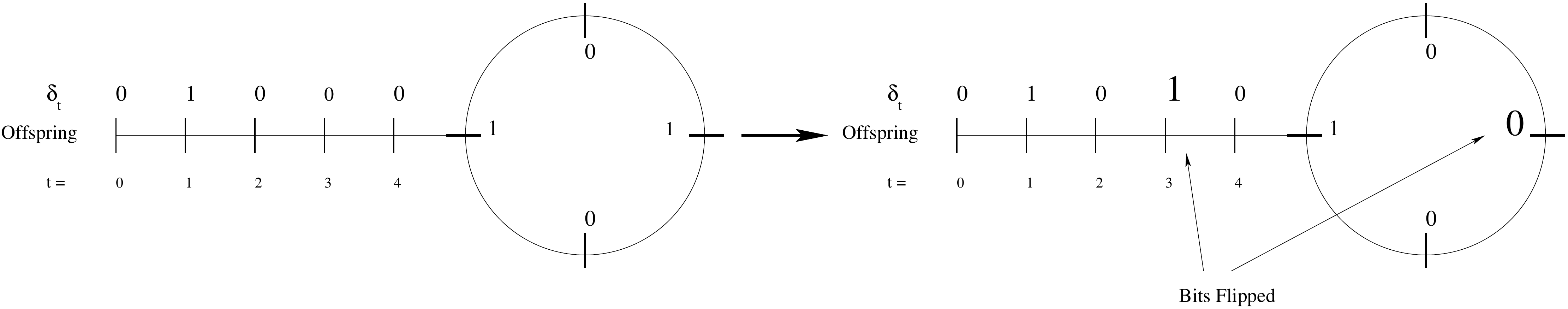,width=0.85\linewidth}
\end{center}
\caption{Type 1 Mutation: Mutation of the offspring strategies by flipping the bits based on mutation probability.}
\label{fig:mutation1}
\end{figure}

\begin{figure}[tbp]
\begin{center}
\epsfig{file=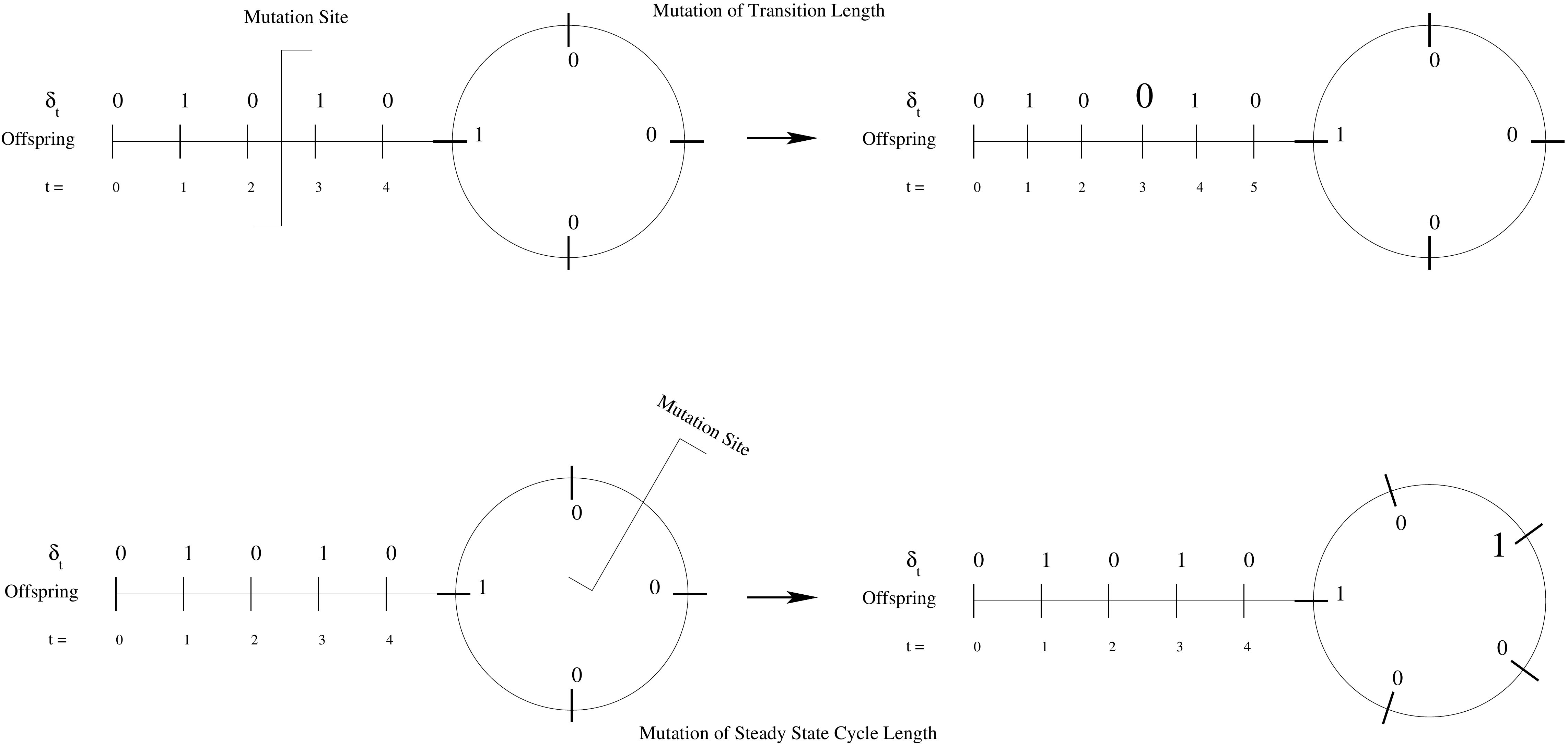,width=0.85\linewidth}
\end{center}
\caption{Type 2 Mutation: Mutation of the transition period length and the steady state cycle length.}
\label{fig:mutation2}
\end{figure}

There are two kinds of mutation that we perform in this study. In type 1 mutation, we go through each genetic bit and flip it from 0 to 1 or 1 to 0 with a probability of $p_m$. This introduces diversity in harvesting strategies and helps in exploring the search space that is not reachable only by sharing of information between the parents. Type 2 mutation is performed to explore an appropriate transition length and cycle length. Type 2 mutation is applied first on the transition string with probability $p_m$. Thereafter, it is applied on the steady state cycle string with probability $p_m$. When this operation is applied, the length of the transition period or steady state period increases or decreases by 1 with equal probabilities. In case the length increases then a bit (0 or 1) is inserted at a random location. If the length decreases then a bit is removed from a randomly chosen location. The type 1 and type 2 mutation operations have been shown in Figure~\ref{fig:mutation1} and~\ref{fig:mutation2} respectively.

\subsubsection{Algorithm Description}
Next, we provide a step-by-step description of the algorithm that is used to explore the harvesting decisions and progress towards the optimum.

\begin{enumerate}
\item [S 1:] Initialize $N$ population members (Refer Section \ref{sec:initialization}) representing different harvesting strategies.
\item [S 2:] Assign fitness (Refer Section \ref{sec:fitnessAssignment}) to each member by solving a non-convex optimization problem at {\em level 2}.
\item [S 3:] Initialize a generation counter: $g \leftarrow 0$.
\item [S 4:] Increment the generation counter by 1: $g \leftarrow g+1$.
\item [S 5:] Choose two parents from the population using tournament 
selection \citep{miller1995}.
\item [S 6:] Perform genetic operations (Refer Section \ref{sec:geneticOperators}) to produce the final offspring.
\item [S 7:] Evaluate each offspring by solving a non-convex optimization problem at {\em level 2}.
\item [S 8:] Choose $\lambda$ random members from the population, and pool them with the offspring.
\item [S 9:] Find the best $\lambda$ members from the pool based on the fitness values and update the population by replacing the chosen members in the previous step with the best $\lambda$ members.
\item [S 10:] If the generation counter ($g$) is less than the maximum number of allowed generations ($g_{\max}$) then go to Step 4; otherwise terminate.
\end{enumerate}

\subsubsection{Parameters}
The optimization task at {\em level 1} requires two kinds of parameters. The first set of parameters are the evolutionary parameters, and the second set of parameters are the problem-specific parameters. The parameter setting used in the experiments performed in this paper are given in Table~\ref{tab:parameters}.

\begin{table}[ht]
\caption{Parameters used during evolutionary optimization}
\label{tab:parameters}
\begin{center}
\begin{tabular}{|l|l|l|l|}
\hline
\multicolumn{2}{|c|}{Genetic Parameters} & \multicolumn{2}{c|}{Problem Parameters} \\ \hline
Population Size ($N$)           & 50       & Min. Transition Length ($t_{\mbox{min}}$)          & 10  \\
Crossover Probability ($p_c$)   & 0.9      & Max. Transition Length ($t_{\mbox{max}}$)          & 25  \\
Mutation Probability ($p_m$)    & 0.1      & Min. Steady State Cycle Length ($s_{\mbox{min}}$)  & 1   \\
Update Parameter ($\lambda$)    & 2        & Max. Steady State Cycle Length ($s_{\mbox{max}}$)  & 10 \\ \hline
\end{tabular}
\end{center}
\end{table}

\subsection{Level 2: Optimization for fitness evaluation}\label{sec:nonConvex}

Our fitness evaluation tasks are non-convex optimization problems which need to be 
solved fast. In this section we define the problem, discuss the sources 
of non-convexities and finally propose an initialization procedure
for optimization with good prospects for the solver to find a global optimum. 

\subsubsection{The problem}

Let $t^0=t_{\mbox{len}}$ denote the end of transition (beginning of first cycle) and $t^1=t_{\mbox{len}}+s_{\mbox{len}}$ the end of first cycle (beginning of second cycle). Given $t^0$, $t^1$ and $\delta_t$, for all $t<t^1$, 
the fitness evaluation problem is to find 
$x_t\geq 0$, for all $t\leq t^1$, 
$h_t\geq 0$ and the auxiliary variables, for all $t<t^1$, to
\be
\max \sum_{0\leq t<t^0} \beta^{t\Delta}c_t +\frac{1}{1-\beta^{(t^1-t^0)\Delta}} \sum_{t^0\leq t<t^1} \beta^{t\Delta}c_t
\label{obj1}
\ee
subject to (\ref{dh})--(\ref{Ct}), (\ref{xt1})--(\ref{Ch}), for $t < t^1$, and
\be
y_{t^0}=y_{t^1}.
\label{steady}
\ee
Here (\ref{steady}) is the steady state condition for state variables. In the objective function (\ref{obj1}), the first term accounts for the transition phase and the second term accounts for the steady phase. The multiplier in front of the summation in the second term accounts
for a geometric series of discounted steady state cash flows.



\subsubsection{Non-convexities}

\begin{figure}[tbp]
\begin{center}
$\begin{array}{ccc}
\includegraphics[height=0.4\linewidth,width=0.4\textwidth,natwidth=610,natheight=642]{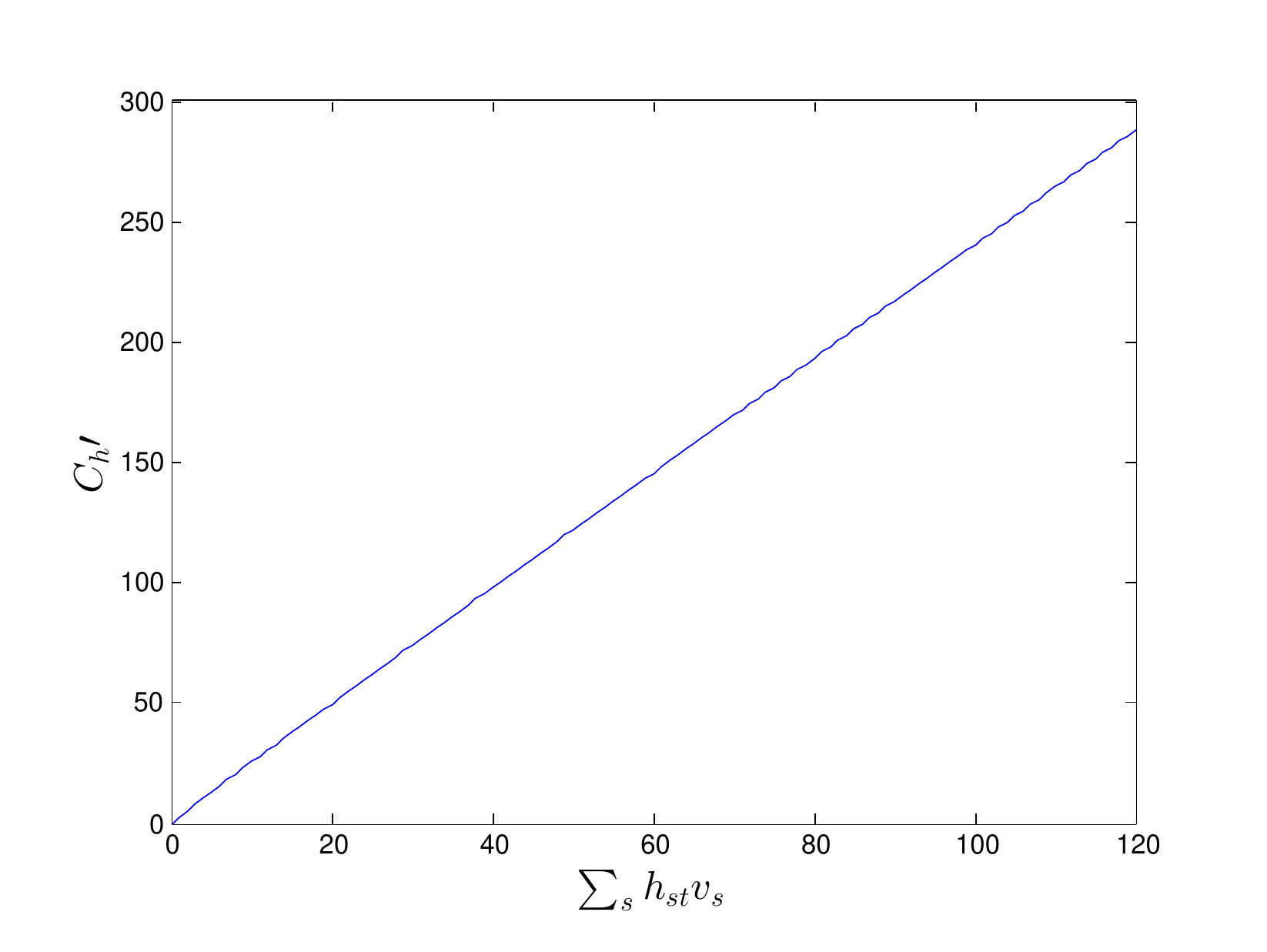}&
\multicolumn{1}{c}{$~~~~$}&
\includegraphics[height=0.4\linewidth,width=0.4\textwidth,natwidth=610,natheight=642]{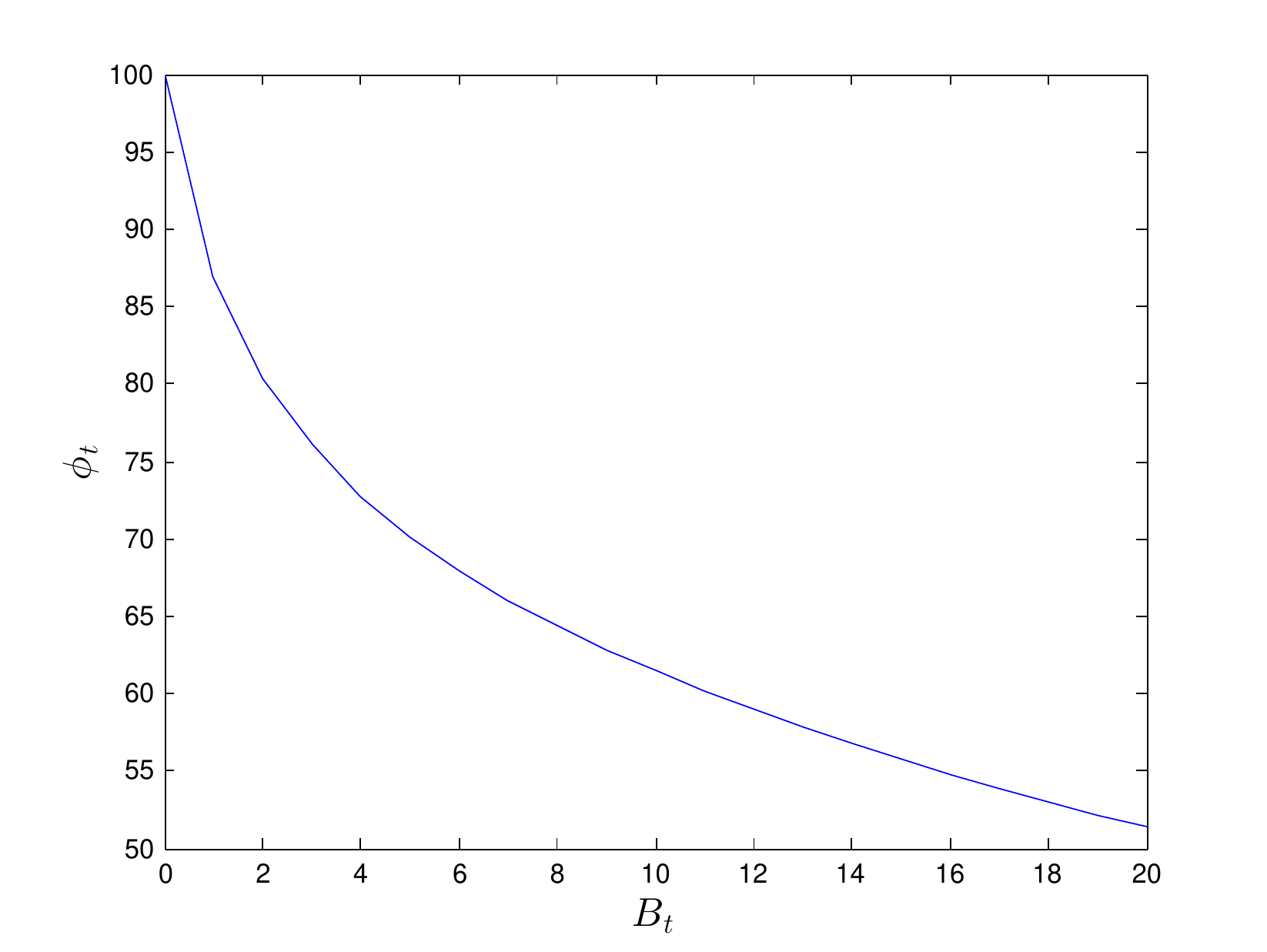}
\end{array}$
\end{center}
\caption{Variable hauling cost as a function of volume $v_t$ in m$^3$ (left) 
and ingrowth as a function of basal area $b_t$ in m$^2$ (right).}
\label{fig:haul_ingr}
\end{figure}

If harvesting takes place at stage $t$, 
then the variable cutting and hauling cost (\euro) in (\ref{Ct}) is $C_c(h_t,v)+C_h(h_t,v)$
where vector $v$ is given in Table~\ref{tab:data}. The graph of
such concave cost function in Figure~\ref{fig:haul_ingr} (left inset) with respect to $h_t$ is almost linear.
In (\ref{st}), ingrowth of trees into size class 1 as a function 
of basal area $B_t$ (m$^2$/ha) is  $\phi_t = (S_1 (B_t+B^0)^{-\nu})/[1+S_2\exp(\gamma B_t)]$ where 
$S_1$, $S_2$, $\gamma$ and $\nu$ are given in Table \ref{tab:data}. 
Figure~\ref{fig:haul_ingr} (right inset) illustrates the ingrowth function for $B_t\in [0,20]$ m$^2$;
for $B_t>10$ m$^2$ the function is mildly non-linear. 

\begin{figure}[!ht]
\begin{center}$
\begin{array}{ccccc}
\includegraphics[height=0.3\linewidth,width=0.3\textwidth,natwidth=610,natheight=642]{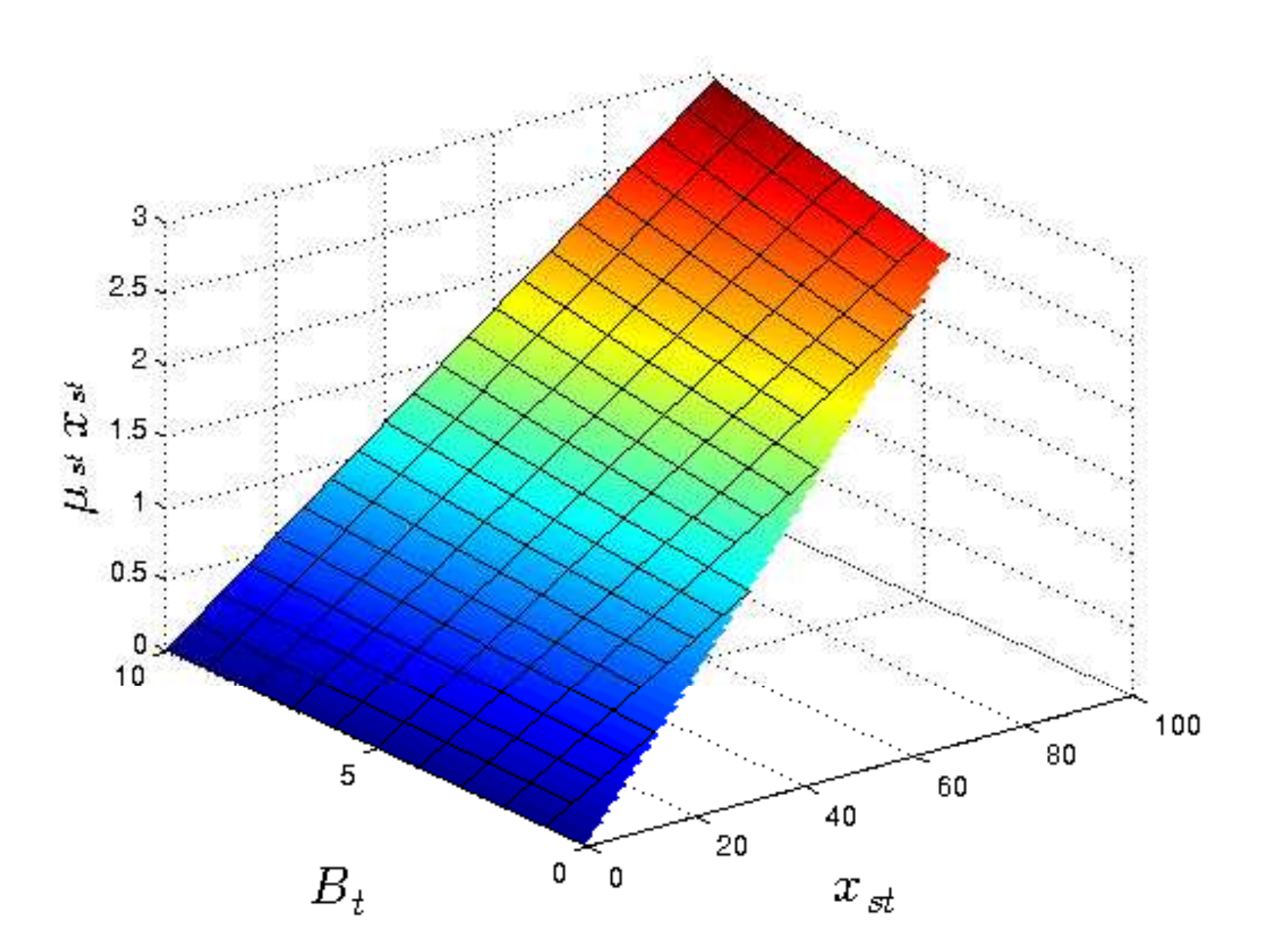}&
\multicolumn{1}{c}{$$}&
\includegraphics[height=0.3\linewidth,width=0.3\textwidth,natwidth=610,natheight=642]{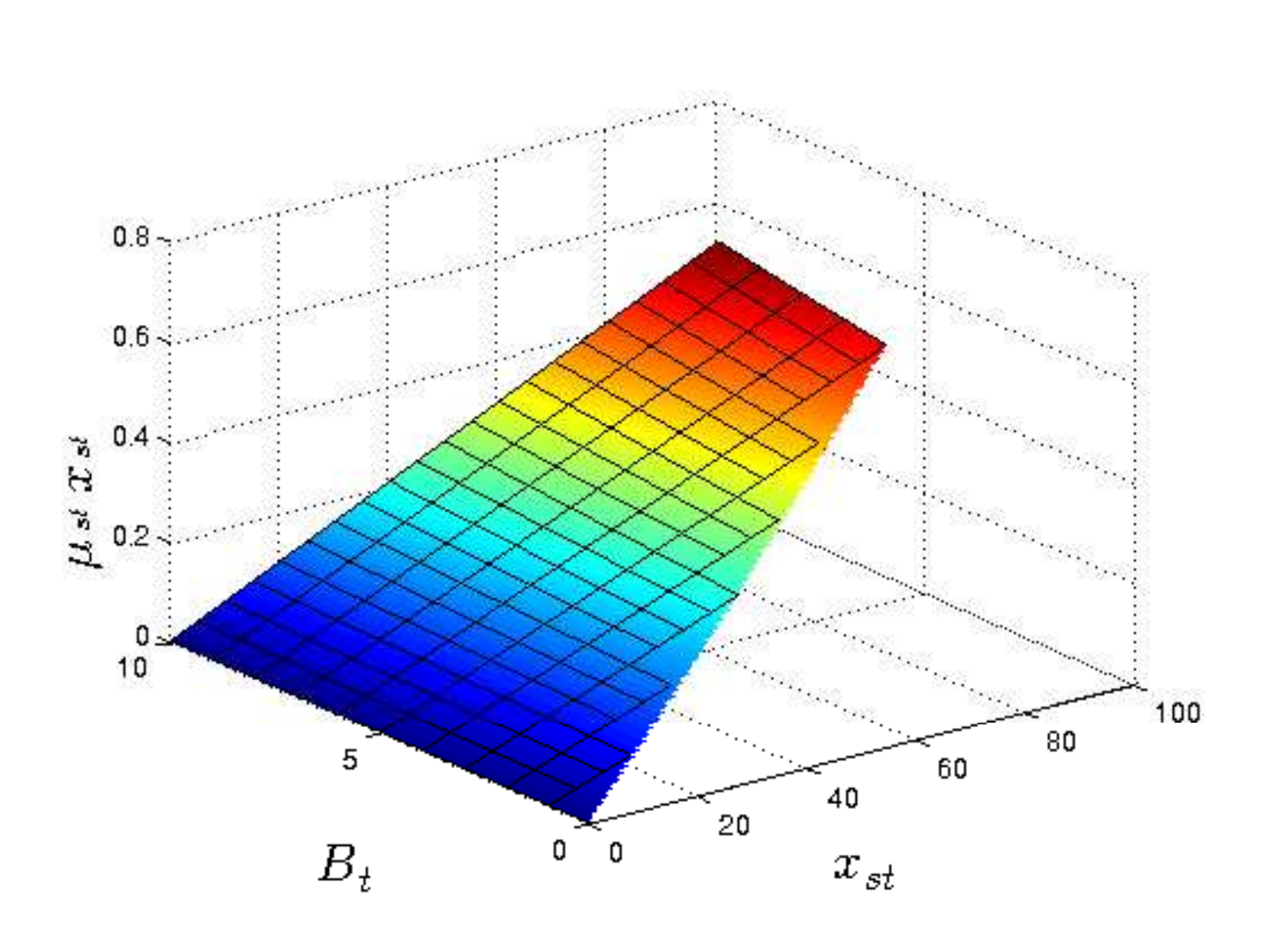}&
\multicolumn{1}{c}{$$}&
\includegraphics[height=0.3\linewidth,width=0.3\textwidth,natwidth=610,natheight=642]{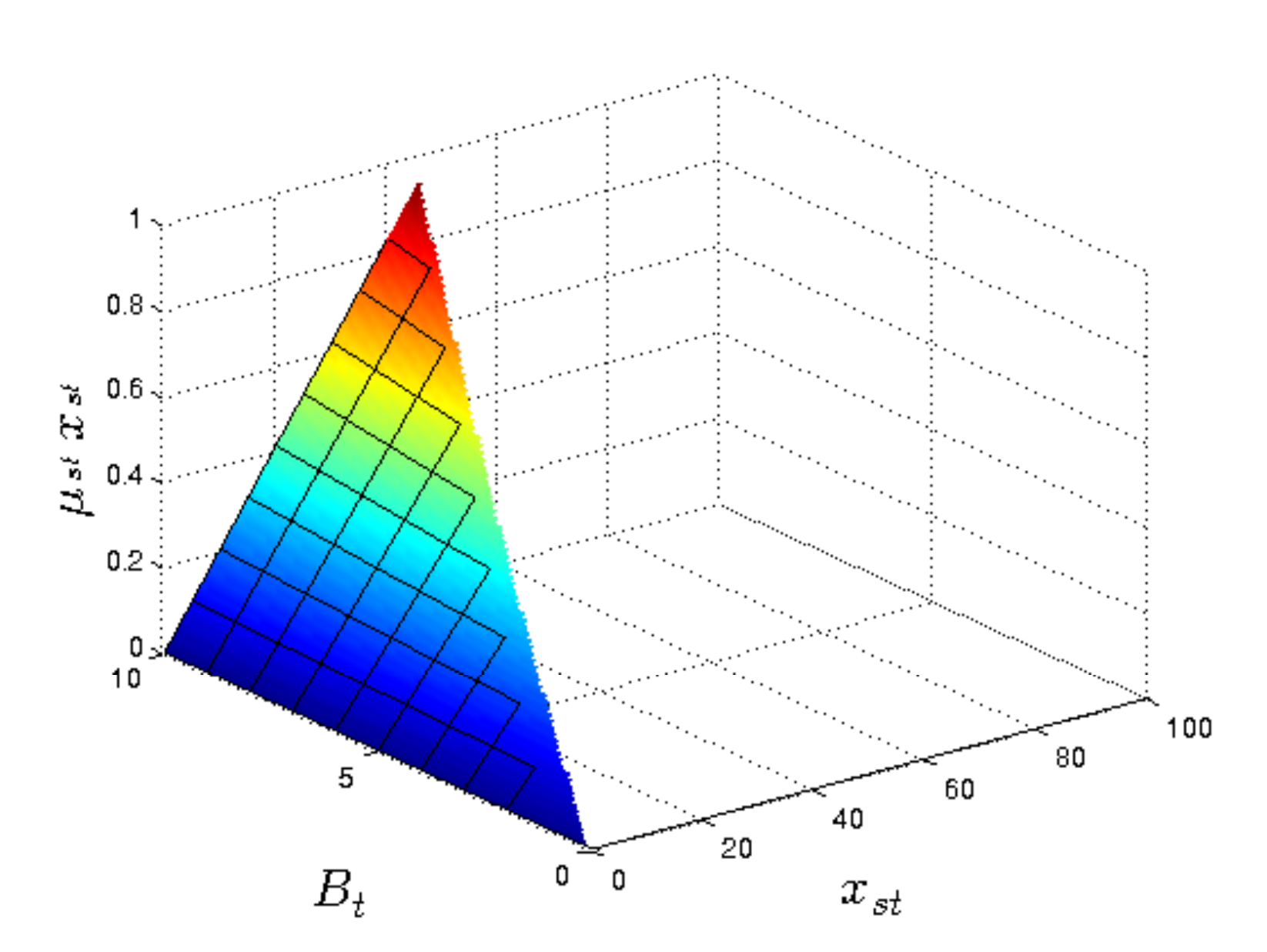}\\
\end{array}$
\end{center}
\caption{Mortality (number of trees) in size classes $s=1$ (left), $s=5$ (middle) and  $s=10$ (right)
as a function of basal area $B_t$ in m$^2$ (left axis) and number of trees $x_{ti}$ (right axis).
}
\label{fig:mor} 

\end{figure}

In (\ref{mti}), mortality fraction of trees dying in size class $s$ is given by 
$\mu_{st} = 1/[1+M_s\exp(-m B_t)]$ where $m$ is given 
in Table \ref{tab:data}. 
Mortality $\mu_{st}$ is a convex function for relevant domain of $B_t$.
In the forest state equations (\ref{xt1})--(\ref{xts}) the reduction 
in the number of trees is $\mu_{st}x_{st}$, where $x_{st}$ is the number of trees in size class $s$.
Figure~\ref{fig:mor} shows such reduction for three size classes $s$, $s=1,5,10$, as a function of $B_t$ 
and $x_{st}$. Note that possible cases are those for which $B_t \geq b_s x_{st}$, 
where $b_s$ is the basal area of a single tree in size class $s$. For such region, 
the mortality functions are almost linear.

\begin{figure}[!ht]
\begin{center}$
\begin{array}{ccc}
\includegraphics[height=0.3\linewidth,width=0.3\textwidth,natwidth=610,natheight=642]{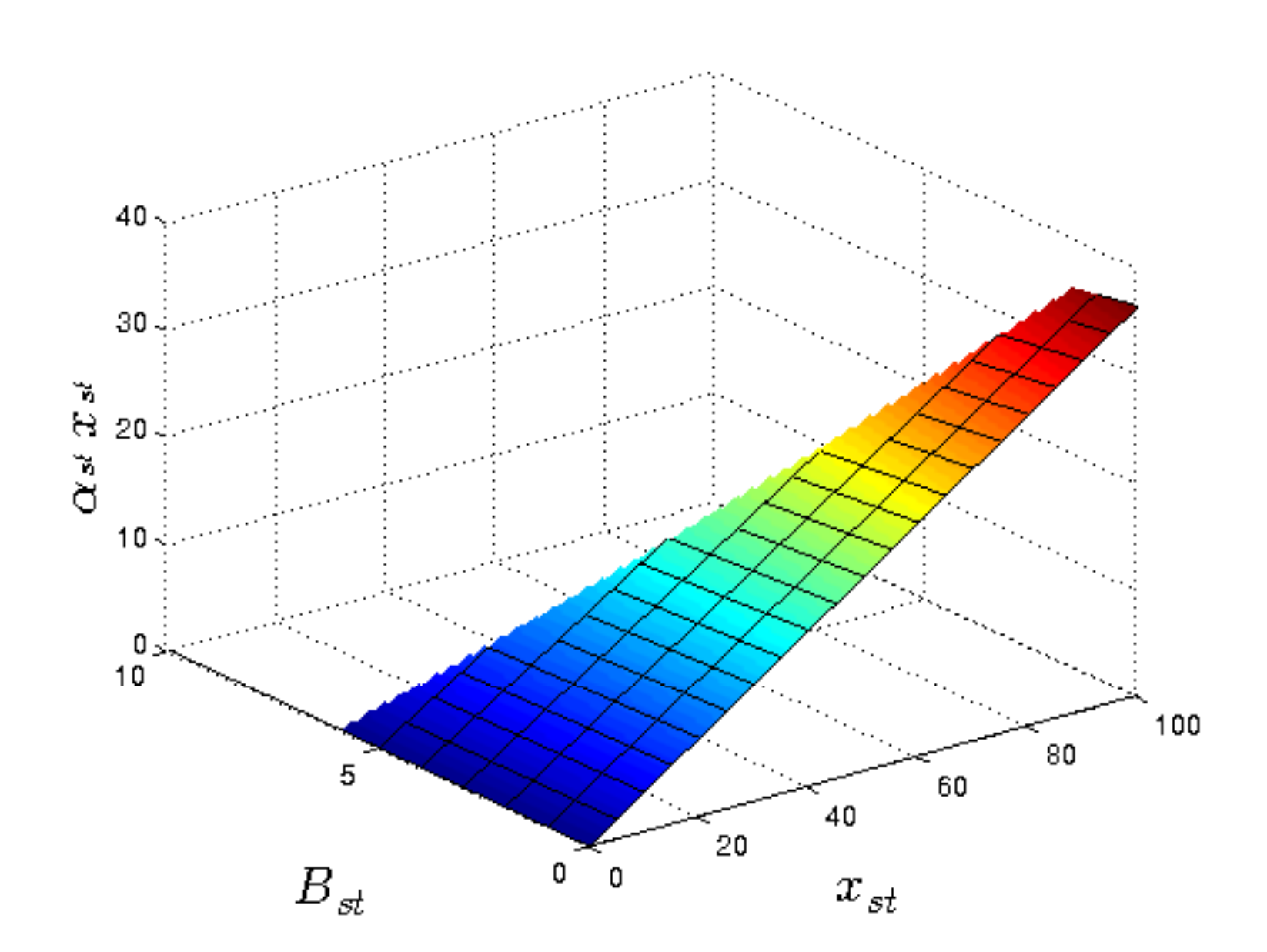}&
\includegraphics[height=0.3\linewidth,width=0.3\textwidth,natwidth=610,natheight=642]{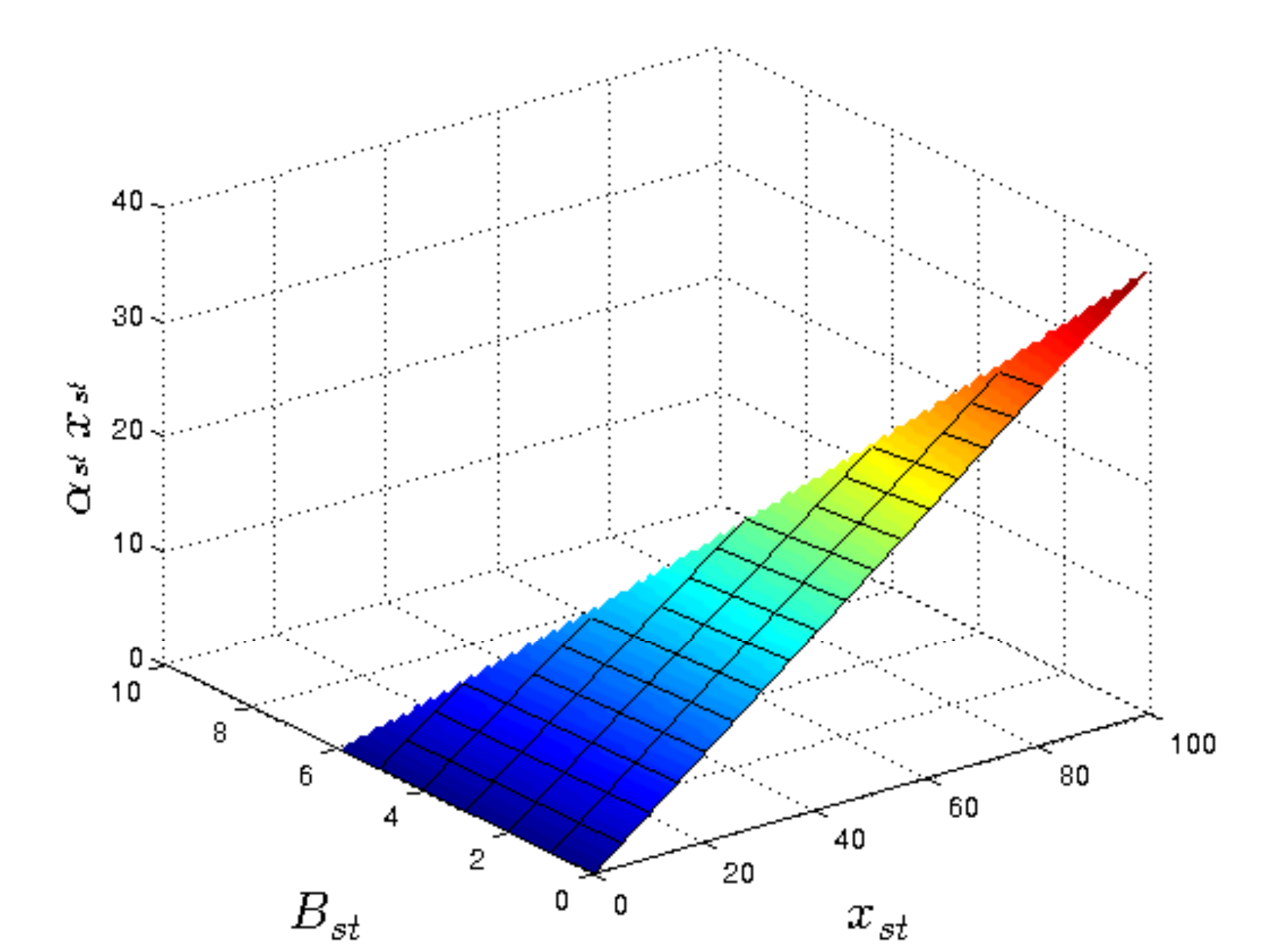}&
\includegraphics[height=0.3\linewidth,width=0.3\textwidth,natwidth=610,natheight=642]{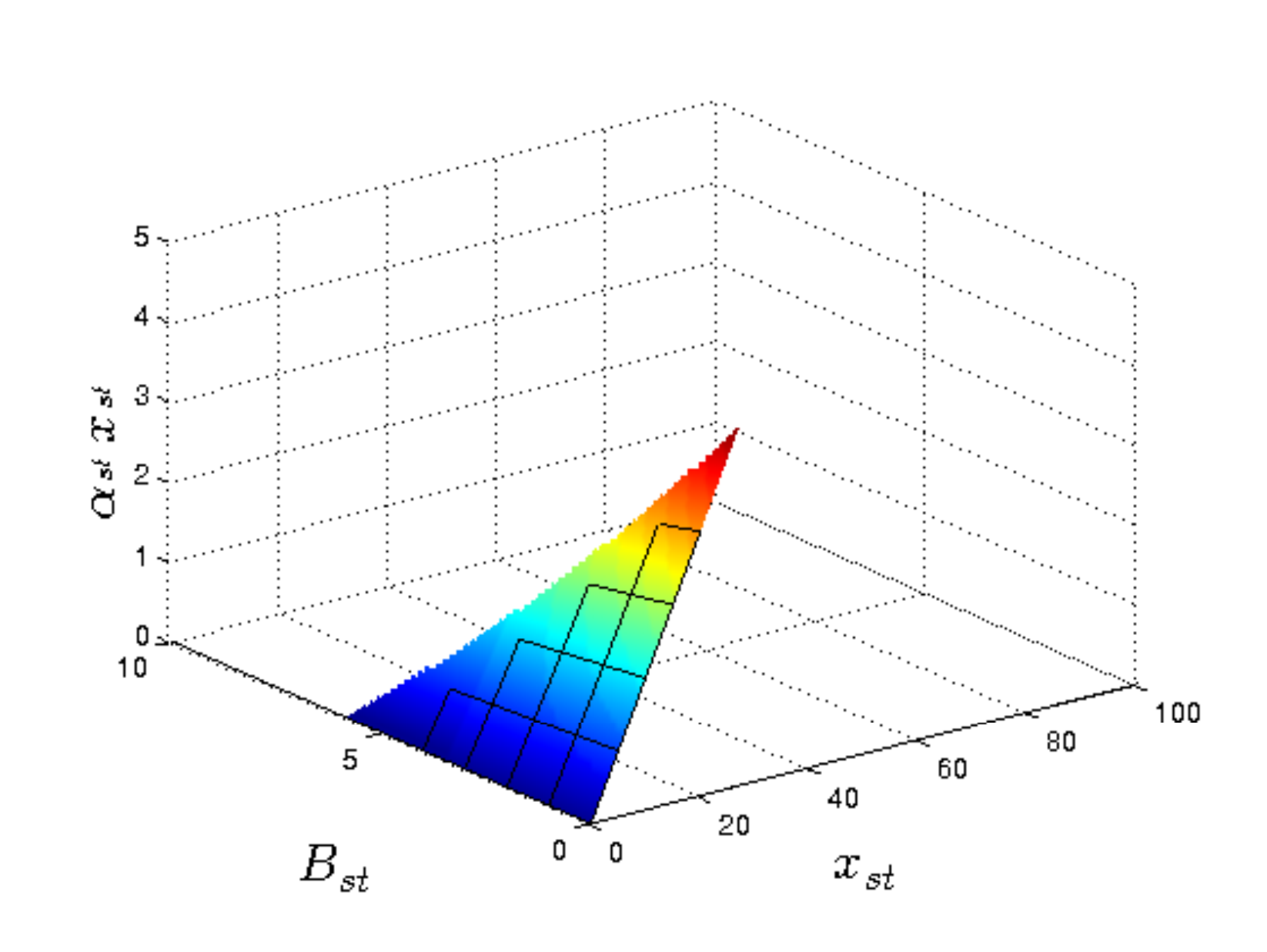}
\end{array}$
\end{center}
\caption{Transition of trees in five years from size class $s$ to $s+1$ 
as a function of the number of trees $x_{st}$
and $B_{st}$.
Cases $s=1$ (left), $s=5$ (middle) and $s=11$ (right) are depicted for basal area $B_t=6$.}
\label{fig:trans}
\end{figure}

For size class $s<12$, given basal area  $B_t$ and 
the basal area $B_{st}$ of trees in size classes above $s$, 
the fraction of trees that grows in five years from size class $s$ to $s+1$ is
$\alpha_{st} = G_s-A_1 B_{st}-A_2 B_t$ where $G_s$, $A_1$ and $A_2$ 
are given in Table \ref{tab:data}. In (\ref{xt1})--(\ref{xts}) 
of the optimization problem the transition 
in terms of number of trees is $\alpha_{st}x_{st}$, which involves a bi-linear function.
Figure~\ref{fig:trans} shows such transitions for three size classes $s$, $s=1,5,10$, 
as a function of $B_{st}$ and $x_{st}$ at $B_t=6$. Possible 
cases are those for which $0\leq B_{st}\leq B_t -b_s x_{st}$ and in such region
again the functions are almost linear.
Given small values of $\rho$ we may term the non-convexities as relatively mild. 
We further confirm this  (refer to Section \ref{sec:randomInit} and Appendix \ref{sec:ec1}) 
by solving the non-convex optimization procedure with different random 
initialization schemes i.e., choosing the starting point for the solver.

\subsubsection{Random initialization for solver}\label{sec:randomInit}

In Appendix \ref{sec:ec1} we test seven alternative versions of random initialization (starting point selection) 
for the solver Knitro \citep{byrd2006knitro} to find a local optimum for the problem (\ref{obj1}).
All seven procedures performed quite well with a high chance for Knitro to end up with a global optimal solution.
Thereby, we chose a randomized starting point for the solver as described in the following. 
Let $s$ denote a size class such that harvesting 
only applies for size classes $i>s$. With $s=5$ exogenously given, let $\eta_{it}$ denote 
the share of trees that are not harvested in size class $i$ at stage $t$. 
Then $\eta_{it}=1$ for $i\leq s$; otherwise  $\eta_{i+1,t}= \epsilon_{it}\eta_{it}$ 
where $\epsilon_{it}$ is drawn from uniform distribution $U(0,1)$. 
Thus the share of trees harvested in size class $i$ at stage $t$ is $(1-\eta_{it})$ 
and it increases in random proportions with $i$, for $i>s$. Shares $\epsilon_{it}$ 
are drawn independently for each time stage $t$ and size class $i>s$.
The number of trees harvested is $h_{it}=\delta_t(1-\eta_{it})x_{it}$ for which we need
the forest state vector $x_t$ unless $\delta_t=0$. Initially, $x_0=x^0$ is given and we obtain 
$h_{i0}=\delta_0(1-\eta_{i0})x_{i0}$
as well as cash flow $C_t$ from (\ref{vt})--(\ref{Ct}). 
Thereafter, auxiliary variables in (\ref{bt})--(\ref{gti}) and 
state equations (\ref{xt1})--(\ref{xts}) yield $x_1$. Similarly, in forward recursion we
obtain values $x_t$, $h_t$, $C_t$ and the auxiliary variables for all $t$. 
In this random initialization procedure we neglect the steady state condition $x_{t^0}=x_{t^1}$.

Let $p$ denote the probability that a random initialization procedure
for the solver ends up with a suboptimal solution. If the problem is initiated independently $k$ times, 
then the probability of not finding the global optimum is $1-p^k$. In our tests reported 
in Appendix \ref{sec:ec1}, $p$ appears to be small (of the order of 0.1 \%). Hence, $1-p^k$ is 
very close to 1 even for $k=1$. This observation may be explained 
by relatively mild non-convexities in our valuation problems; 
see nonlinear function illustrations in Figures~\ref{fig:haul_ingr}-\ref{fig:trans}.

\bigskip


\section{Results}\label{sec:results}
In this section, the proposed algorithm has been applied on an empirically estimated size-structured ecological model for uneven-aged Norway spruce forests \citep{bollandsas2008}. We have considered three different initial states of the forest, as stated in Table \ref{tab:data_i}. The first initial state ($x_0=x^1$) is a young even-aged stand, the second ($x_0=x^2$) is an uneven-aged stand, and the third ($x_0=x^3$) is an old even-aged stand. All the forest parameters and functions in our study are fixed as given in Tables \ref{tab:data_i} and \ref{tab:data}. Sensitivity studies have been performed by varying the parameters: interest rate ($r$), fixed cost ($C_f$) and site index ($S$). The different values of these parameters considered are given in Table \ref{tab:sensitivity}.

\begin{table}[h]
\begin{center}
\label{tab:sensitivity}
\caption{Different values of the parameters studied in the paper corresponding to 3 different initial states of the forest}
\begin{tabular}{lccc}
\toprule
            & $r$ ($\%$) & $C_f$ (\euro) & $S$  \\
\midrule
Cases studied &  1,2,3,4   &  100,300,500     &  13,15,17   \\
Base Case   &  3   &   300    &  15    \\
\bottomrule
\end{tabular}
\end{center}
\end{table}

It was observed from our runs that the optimal steady state harvesting interval and size distribution is independent of the initial stand state of the forest. Steady state results corresponding to the three initial states and different values of parameters $r$ and $C_f$ are shown in Table \ref{tab:resultsTable}. Site index, S, was kept fixed at 15 in this study. Some of the important observations drawn from these runs have been outlined below:
\begin{enumerate}
\item The optimal solution converges to the steady state harvesting interval within 120 years (24 intervals) or sooner (Table \ref{tab:resultsTable}). 
\item The transition towards the steady state tree size distribution is shortest when the initial stand state is already heterogeneous and longest for the young even-aged stand. 
\item The length of steady state interval varies between 10 and 25 years and increases or remains the same for higher levels of fixed cost.
\item Within the Faustmann optimal rotation framework higher interest rate implies shorter rotation suggesting that the length of the steady state interval decreases in the interest rate. However, this is not the case: for example, when fixed cost is equal to \euro 100 and interest rate is increased from $3\%$ to $4\%$ harvesting interval lengthens from 10 to 15 years. 
\item When interest rate increases it is optimal to allocate a larger fraction of capital from forestry to alternative sources. This is reflected in the fact that when harvesting interval lengthens with interest rate, it becomes optimal to cut smaller size class (i.e., size class with average diameter equal to 225mm) and the stand volume both before and after harvest decreases. 
\item Given any fixed cost level the stand volume before and after harvest is lower or the same with higher interest rate.
\end{enumerate}

\begin{sidewaystable}
\begin{center}
\caption{Steady State (s-s) Results for different values of $r$ and $C_f$ with $S$ fixed at 15}
\label{tab:resultsTable}
\begin{tabular}{ccccccccc}
\toprule
$C_f$	&	$r$	&	Interval between	&	$\frac{\mbox{Profit}}{\mbox{Year}}$	&	$\frac{\mbox{Volume}}{\mbox{Harvest}}$	&	$\frac{\mbox{Avg. Volume}}{\mbox{Year}}$	&	Harvested	&	Period when s-s & No. of trees	\\

	&		&	harvests in s-s	&		&		&		&	threes	&	interval reached for & before and\\

 (\euro)	&		&	(years)	&	(\euro)	&	$(m^3)$	&	$(m^3)$ & Size (mm) 	&	initial states $x^1, x^2, x^3$	& after harvest\\
\midrule
100	&	0.01	&	15	&	308	&	91.6	&	6.1	&	325-425	&	18,4,24	&	824/715	\\
100	&	0.02	&	15	&	265	&	81.3	&	5.4	&	324-425	&	9,0,16	&	753/618	\\
100	&	0.03	&	10	&	250	&	52.4	&	5.24	&	275-325	&	7,0,18	&	715/620	\\
100	&	0.04	&	15	&	193	&	63.8	&	4.3	&	225-325	&	6,0,12	&	667/505	\\
300	&	0.01	&	15	&	295	&	91.6	&	6.1	&	325-425	&	19,2,18	&	824/715	\\
300	&	0.02	&	20	&	266	&	111.3	&	5.6	&	275-425	&	9,1,12	&	788/616	\\
300	&	0.03	&	15	&	252	&	81.3	&	5.4	&	275-375	&	13,0,18	&	753/618	\\
300	&	0.04	&	20	&	201	&	90.2	&	4.51	&	225-375	&	7,0,18	&	709/502	\\
500	&	0.01	&	20	&	291	&	121.7	&	6.1	&	325-475	&	12,3,10	&	849/710	\\
500	&	0.02	&	20	&	256	&	111.3	&	5.6	&	275-425	&	9,1,17	&	788/616	\\
500	&	0.03	&	25	&	210	&	118.3	&	4.7	&	225-425	&	9,0,17	&	748/499	\\
500	&	0.04	&	20	&	191	&	90.2	&	4.51	&	225-375	&	8,0,12	&	709/502	\\
\bottomrule
\end{tabular}
\end{center}
\end{sidewaystable}

\subsection{Dependence of Optimal Solution on Initial State}
\begin{figure}[!ht]
\begin{center}
\epsfig{file=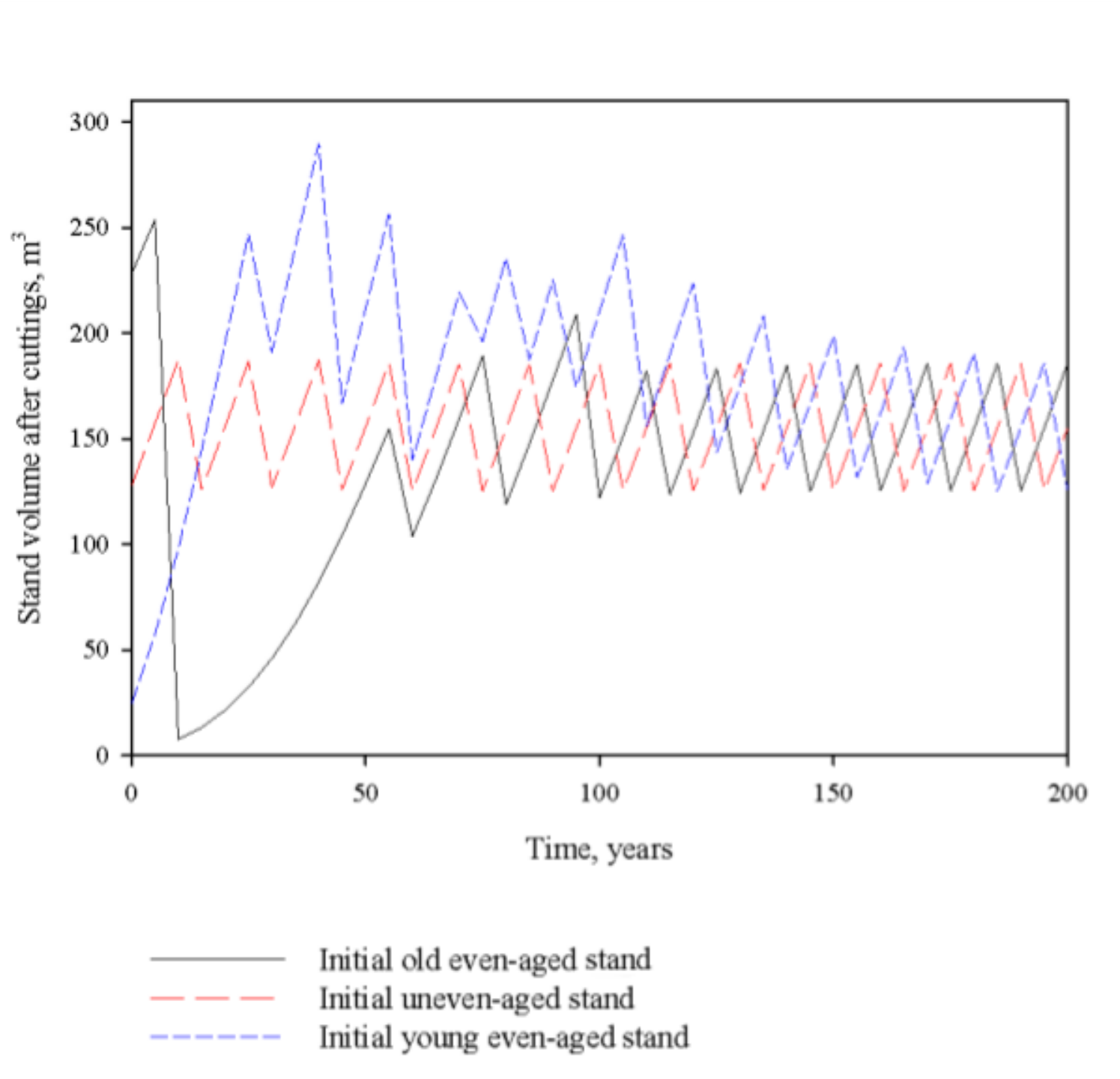,width=0.6\linewidth}
\end{center}
\caption{Dependence on the initial state ($r=0.01, C_f=300$).}
\label{fig:resultsFig1}
\end{figure}

The optimal solution depends strongly on the initial stand state (Figure \ref{fig:resultsFig1}). Given an initially dense stand (solid line) it is optimal to almost clearcut after 10 years and then wait 50 years until the next harvest. In contrast when the initial stand is already heterogeneous (Table \ref{tab:data_i}) the steady state harvesting interval (15 years) is optimal immediately from the beginning even if it takes about 140 years to reach the exact steady state tree size distribution and harvest level. Given an initial stand containing trees only in the smallest size class, it is optimal to wait 25 years until the first harvest and it takes 80 years to reach the steady state harvesting interval.

\subsection{Dependence of Optimal Solution on Interest Rate and Fixed Costs}
\begin{figure}[!ht]
\begin{center}
\epsfig{file=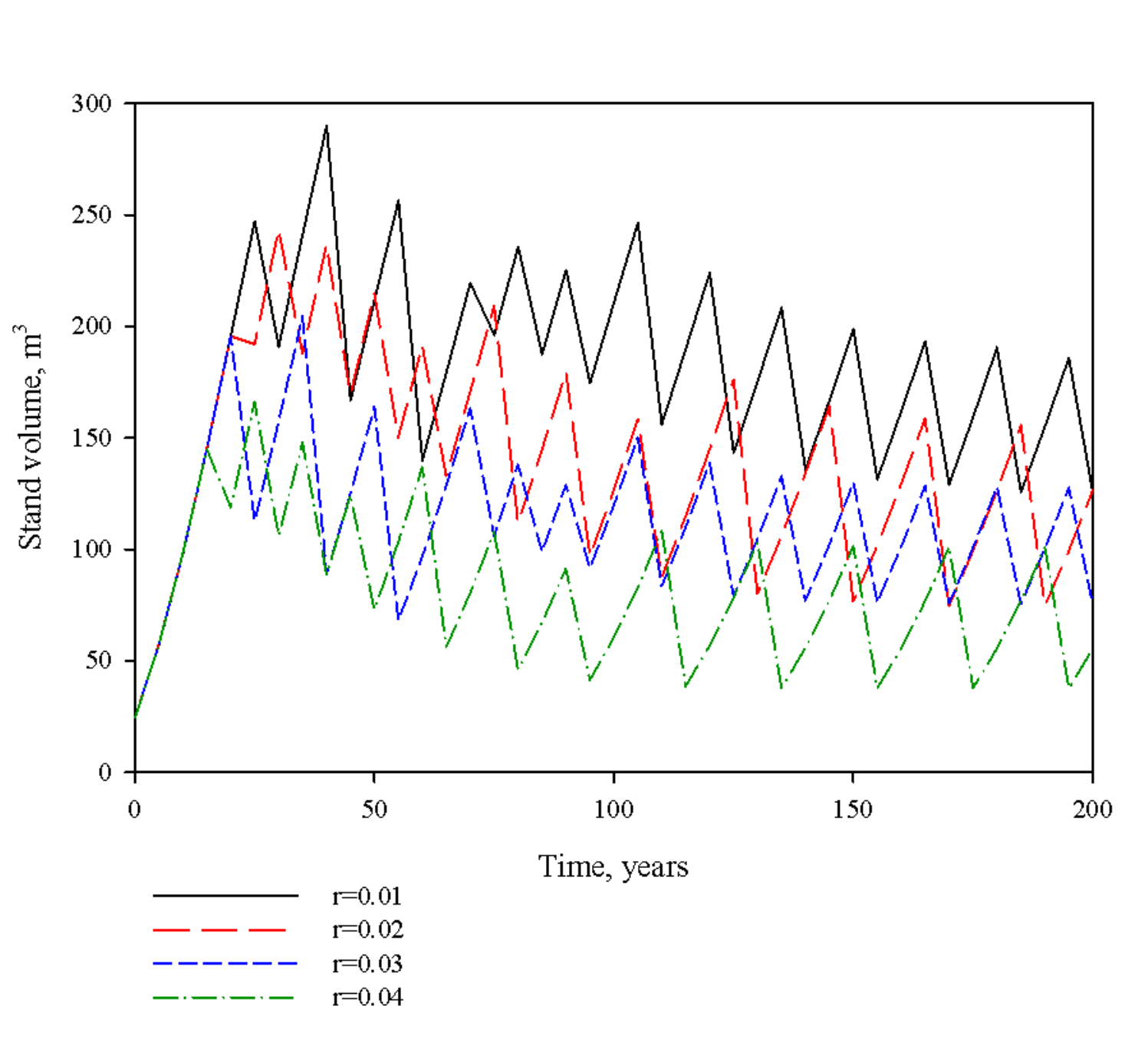,width=0.6\linewidth}
\end{center}
\caption{Dependence of optimal solution on interest rate (Initial state $x_0 = x^1$, $C_f=300$).}
\label{fig:resultsFig2}
\end{figure}

\begin{figure}[!ht]
\begin{center}
\epsfig{file=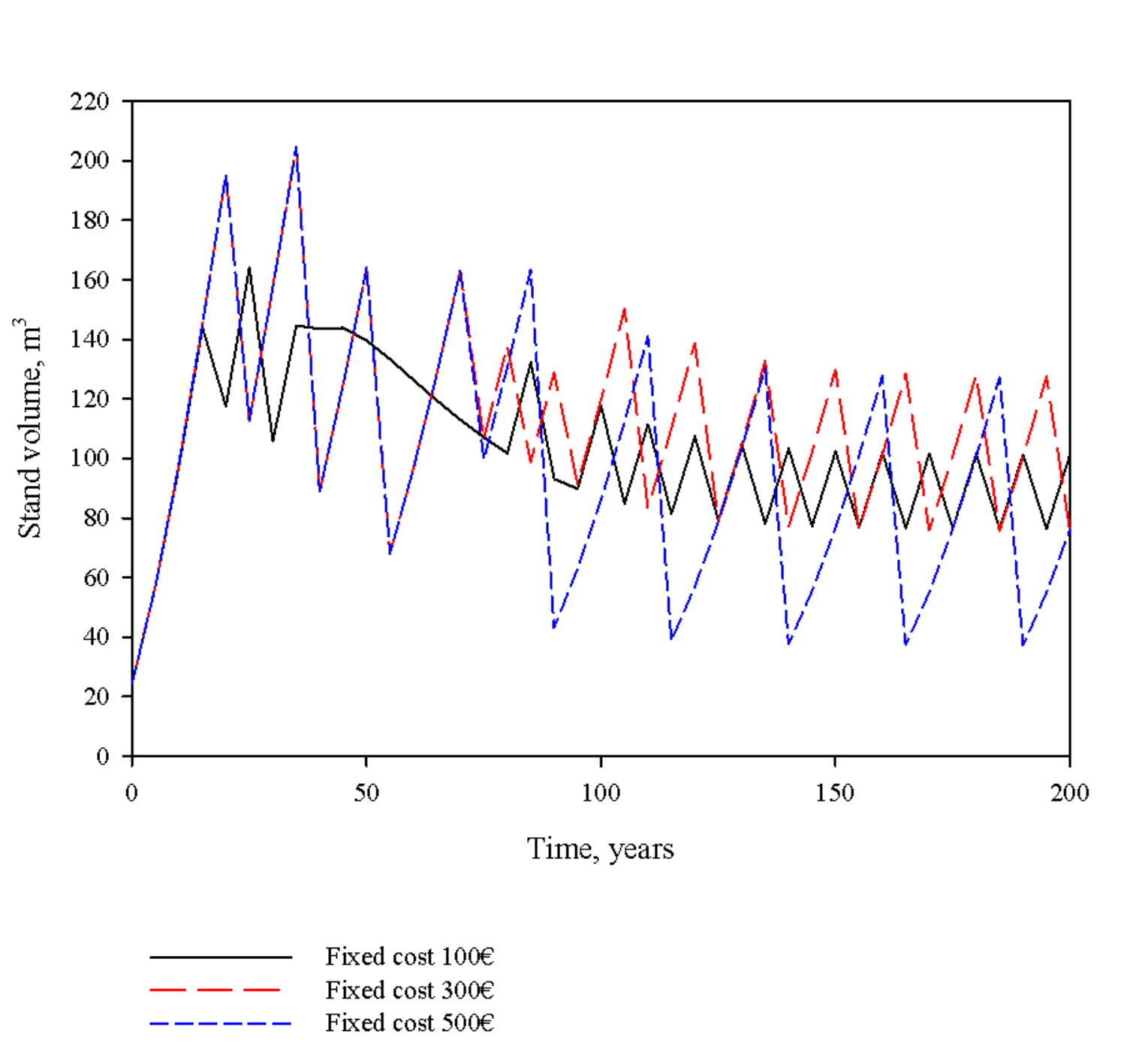,width=0.6\linewidth}
\end{center}
\caption{Dependence of optimal solution on fixed cost (Initial state $x_0 = x^2$, $r=0.03$).}
\label{fig:resultsFig3}
\end{figure}

In Figure \ref{fig:resultsFig2} the initial stand contains only young trees and 
the first cutting is postponed to 30, 25 and 20 years when the interest rate is 
increased from $1$ to $4\%$. Varying the fixed cost level above \euro 100 has 
strong effects on optimal harvest timing (Figure \ref{fig:resultsFig3}). 
Given an initially young stand, interest rate equal to $3\%$ and \euro 100  
fixed cost, there is a 45 year period with harvest every 5 years, while with
 fixed cost equal to \euro 500  it is always optimal to wait at least 15 years 
before the next cut. Fixed cost must be decreased to \euro 20  until the steady 
state solution is to harvest the stand every period.

\subsection{Dependence of Optimal Solution on Site Index}
\begin{figure}[!ht]
\begin{center}
\epsfig{file=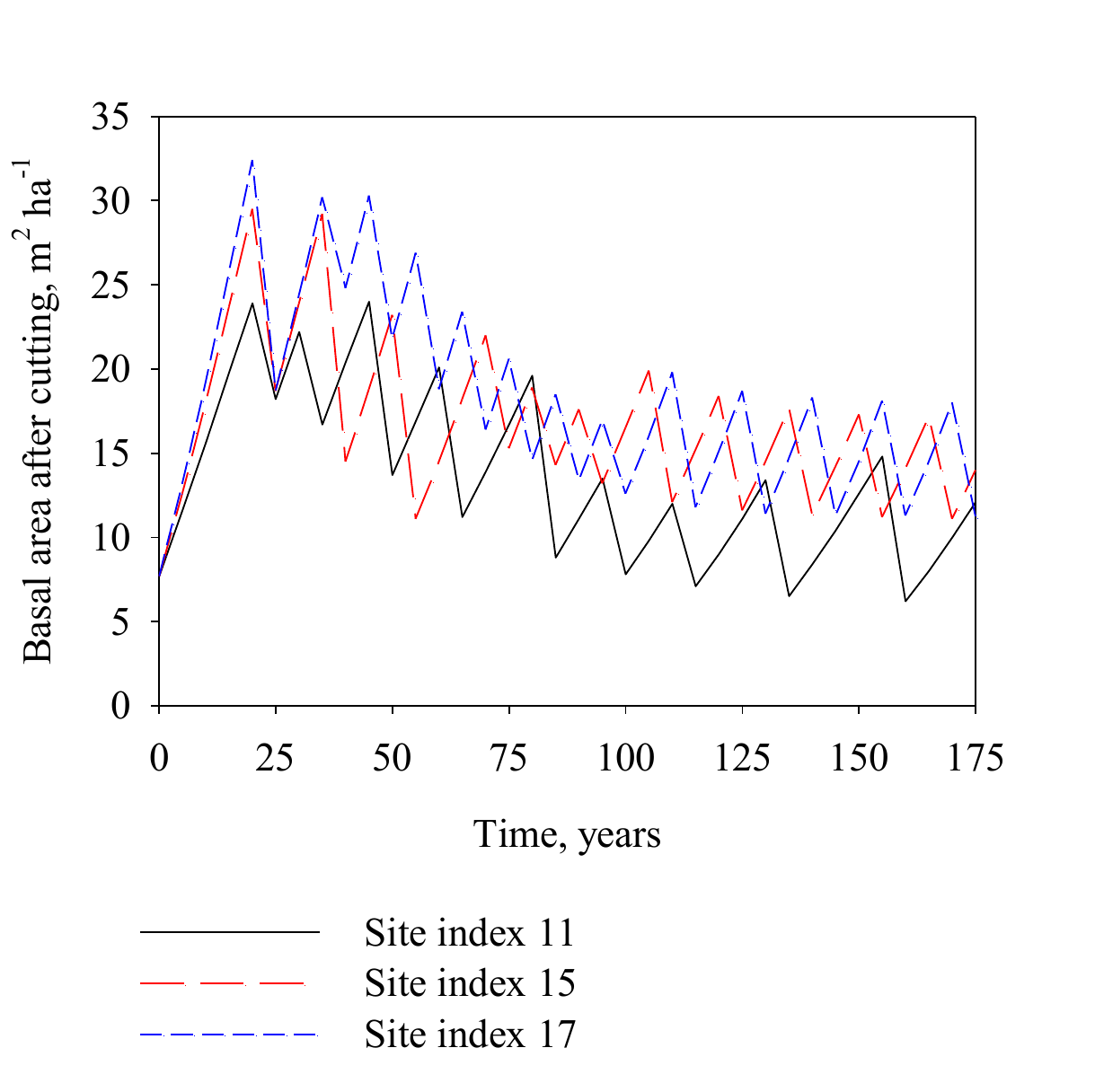,width=0.6\linewidth}
\end{center}
\caption{Dependence of optimal solution on the site index (Initial state $x_0=x^1$, $C_f=300$, $r=0.03$)}
\label{fig:siteIndex}
\end{figure}

Site index has as strong effect on the the optimal transitioning and steady state that can be observed in Figure \ref{fig:siteIndex}. Given that initial stand contains trees only in the smallest size class, the time of first harvest is at year 20, regardless of site index. The stand density at the first harvest; however, differs between the sites; with more productive sites the overall density of the stand is higher compared to less productive sites. Similarly with other initial stand states, the timing of the first harvest is the same regardless of site index. However, the timing of the second and subsequent harvests during the transition as well as in the steady state differ between the site types. At less productive sites the harvests occur less frequently and target smaller size classes compared to more productive sites. During transition, the more productive sites are kept at higher densities, but in the steady state both the harvest timing and size of harvested trees are the same with site indices 15 and 17. However, annual profit and harvested volume are higher at the more productive site due to the higher growth rate (refer to Table \ref{tab:siteIndex}).

\begin{table}[h]
\begin{center}
\footnotesize
\caption{Steady state results for different site indices (Initial state $x_0=x^1$, $C_f=300$, $r=0.03$)}
\label{tab:siteIndex}
\begin{tabular}{ccccccc}
\toprule
Site	&	Interval	&	$\frac{\mbox{Profit}}{\mbox{Year}}$	&	$\frac{\mbox{Volume}}{\mbox{Harvest}}$	&	$\frac{\mbox{Avg. Volume}}{\mbox{Year}}$	&	Harvested trees	&	 No. of trees before	\\
Index		&	(years)	&	(\euro)	&	$(m^3)$	&	$(m^3)$ & Size (mm) 	&	 and after harvest\\
\midrule
11	& 25	& 138	& 78	& 3.1	& 225-425	& 694/498 \\ 
15	& 15	& 252	& 81	& 5.4	& 275-375	& 753/618 \\
17	& 15	& 289	& 92	& 6.1	& 275-375	& 770/620 \\
\bottomrule
\end{tabular}
\end{center}
\end{table}

\subsection{Other Results}
\begin{figure}[!ht]
\begin{center}
\epsfig{file=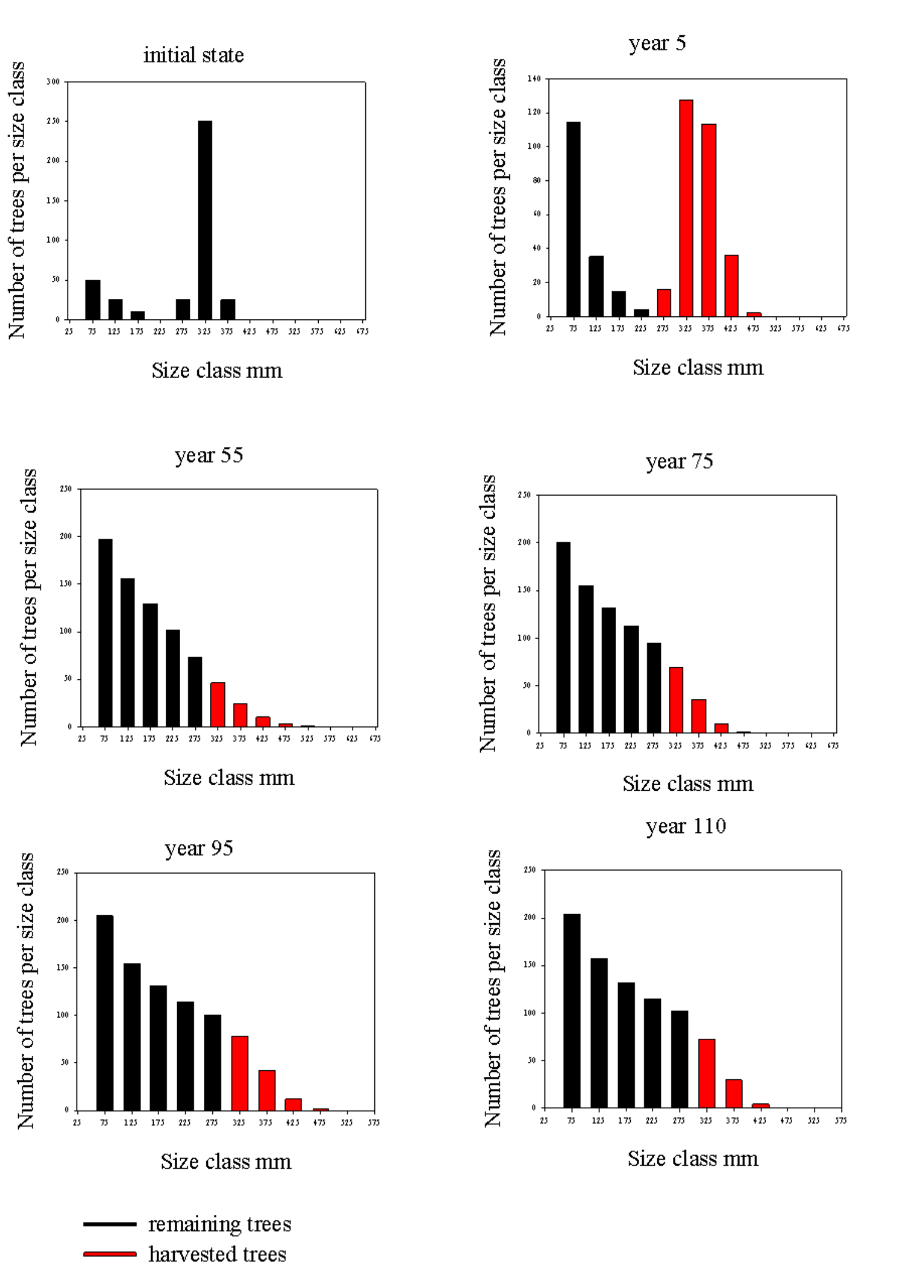,width=0.85\linewidth}
\end{center}
\caption{Development of size distribution over time (Initial state $x_0 = x^3$, $r=0.01$, $C_f=300$).}
\label{fig:resultsFig4}
\end{figure}

Figure \ref{fig:resultsFig4} shows the results for an initially uneven-aged stand ($x_0=x^2$). 
We observe that the stand develops close to the steady state values in 95 years. 
When steady state harvesting interval is three periods, the steady state harvest is 
targeted to three oldest size classes. However, during the transition phase with 
varying harvesting interval, trees are harvested from size classes 3-5. 
\color{black}Economic optimization tends to yield solutions where harvesting applies
to larger trees and no smaller trees are harvested when natural regeneration is 
insufficient.\color{black}

The economic gain from optimizing the harvesting interval against a fixed harvesting interval was also compared for different cases. The gain was found to be remarkable when the initial state is far from the steady state. For cases shown in Figure \ref{fig:resultsFig2}, an optimal harvesting schedule increases the objective function value by about $10-12\%$.

All the results presented in this section so far are based on the computations performed using our proposed algorithm. The computational efficiency of the approach allowed us to perform large number of runs for a variety of cases, that would have been difficult otherwise. In order to gain confidence in the results achieved by our algorithm, we evaluate  our approach against the standard branch-and-bound algorithm that is commonly used to handle mixed integer programs.

\subsubsection{Comparison of Evolutionary Optimization with Branch-and-Bound}\label{sec:comparisons}

\color{black}In Appendix B we demonstrate how the problem can be solved exactly using mixed integer programing with branch-and-bound.  In Table 8 we provide a comparative evaluation of our heuristic evolutionary optimization algorithm (EO) against the exact   branch-and-bound approach (BB) for the base case of 3\% interest rate ($r$) and \euro$~$300 fixed cost $C_f$ . Comparisons have been drawn in terms of the best objective value achieved (1000\euro), run time (hours) and number of non-convex         optimization calls made. The evolutionary optimization algorithm was terminated when a maximum number of non-convex optimization calls exceeded 8000, while the branch-and-bound approach was terminated either when the maximum number of      active nodes in the method exceeded 100,000 or the method exited normally. In Table 8, termination of EO due to maximum number of non-convex optimization calls is denoted as NCO, termination of BB due to maximum number of active nodes is   denoted as NAN, and termination of BB at confirmed optimum is denoted as NOR. For further information about the mixed integer programing solved using branch-and-bound, the readers may refer to Appendix B.

        The results in Table 8 suggest that evolutionary optimization approach offers an enormous computational advantage over the branch-and-bound approach. In terms of the best objective value achieved, evolutionary optimization algorithm was found to   outperform branch-and-bound marginally for initial state $x_2$. In this case BB was terminated by the upper limit on the number of active nodes. For initial states, $x_1$ and $x_3$, the difference in the best NPV found by EO and BB is truly insignificant.  The run time results are not surprising given that branch-and-bound performs an exhaustive search in the decision space.\color{black}


\begin{table}[h]
\caption{Comparison between evolutionary optimization (EO) approach proposed in this paper against the branch-and-bound (BB) approach for the base case of $3\%$ interest rate and \euro 300 fixed cost. Objective values are in 1000\euro/ha}
\begin{center}
\label{tab:comparisonTable}
\begin{tabular}{|c|c|c|c|c|c|}
\hline
\begin{tabular}[c]{@{}c@{}}Initial\\ Condition\end{tabular} 
& Method 
& \begin{tabular}[c]{@{}c@{}}Best\\ Objective\\ Achieved\end{tabular} 
& \begin{tabular}[c]{@{}c@{}}Time\\ Required\\ (Hours)\end{tabular} 
& \begin{tabular}[c]{@{}c@{}}Non-Convex\\ Optimization\\ Calls\end{tabular} 
& Termination \\ \hline
\multirow{2}{*}{$x^1$} & EO & 10.5535 & 2.14 & 8,000 & NCO \\ 
& BB & 10.5534 & 41.42 & 137,595 & NAN \\ \hline
\multirow{2}{*}{$x^2$} & EO & 8.9632 & 2.21 & 8,000 & NCO \\ 
& BB & 8.9448 & 52.61 & 165,047 & NAN \\ \hline
\multirow{2}{*}{$x^3$} & EO & 15.1595 & 2.25 & 8,000 & NCO \\ 
& BB & 15.1604 & 86.82 & 258,650 & NOR\\ \hline
\end{tabular}
\end{center}
\end{table}

\section{Conclusions and Future Work}\label{sec:conclusions}
In this paper, we make a two-fold contribution to the uneven-aged forest management literature. 
The first contribution is on the method side where we propose a customized algorithm to handle mixed integer non-linear optimization problem on optimal management of naturally regenerating uneven-aged forests. 
\color{black}For comparison with our heuristic approach we use branch-and-bound, which is 
a method for which optimality can be proven. Our approach is able to 
computationally surpass branch-and-bound by more than an order of magnitude,     
which would facilitate solving large scale forest management problems. 
A fast algorithm enabled us to make a more significant second 
contribution that is building a better insight for the uneven-aged forest model.
A thorough analysis of an empirical model for uneven-aged spruce forest stands
provided the following important results contributing towards uneven-aged management literature: 
\begin{itemize}
\item Optimizing the harvesting interval yields up to $12\%$ increase in net present value income compared to a solution where the harvesting interval is fixed. The efficiency gain is highest when interest rate and fixed cost is high and initial stand is young.
\item Fixed cost was found to have major implications on the optimal solution. For a given initial size class distribution, and interest rate, varying the fixed costs can lead to very different optimal solutions.
\item The stationary state of the stand was found to be independent of the initial stand state.
\item Stationary harvesting interval tends to be longer with higher fixed cost but may lengthen or shorten with interest rate.
\item Stationary harvesting interval is reached much sooner than stationary size class structure.
\item A parametric study for the site index that denotes land fertility suggests that with higher yields the harvesting interval is shorter.
\end{itemize}
The result on efficiency gain when compared with fixed interval harvesting could 
be an important factor for making a shift towards uneven-aged forest management practices, 
given that uneven-aged management also offers the advantage of biodiversity preservation. 
Other results on interest rate study, fixed cost study and site index study 
provided in the paper would be useful in directing the future research on uneven-aged 
forest management. The proposed method would be directly applicable for multi-species 
forests as well, and we intend to study optimal management of such forests in our future work. 
\color{black}

\section{Acknowledgments}
\color{black}The authors are grateful to two anonymous reviewers for constructive comments.~\color{black}
Ankur Sinha and Pekka Malo would like to acknowledge the support provided from Liikesivistysrahasto and Helsinki School of Economics Foundation.

\bibliographystyle{elsarticle-harv}
\bibliography{forest}



\newpage
\begin{appendices}


{\bf Computational tests for fitness evaluation and formulation for mixed integer programing}
\\
In the appendices we discuss tests with alternative random initialization 
procedures for the fitness evaluation problem (Equation \ref{obj1}) and provide mixed integer programing formulation for the problem that has been solved using branch-and-bound algorithm (Section \ref{sec:comparisons}).

\section{Testing random initialization procedures}\label{sec:ec1}

For tests we considered seven versions of random initialization procedures.
Each one has some charter for solving the problem such as harvesting 
applies primarily to mature trees (high size classes). For versions other than 
the initialization explained in Section~\ref{sec:nonConvex}, 
the smallest size class subject to harvesting is varied, 
for instance, it may be randomly drawn. Also the strategy
for the shares of trees harvested varies: e.g., harvesting applies to largest sizes
classes until a target level of  harvested volume is met. Such target level 
may be tied to the volume growth of the forest adjusted by a random perturbation.
In our tests all such initialization procedures performed well in the sense that 
suboptimal solutions were infrequently produced by the solver.


Initial forest states  and harvesting schedules used for tests are illustrated 
in Tables~\ref{tab:yinit}-\ref{tab:harv}
Table~\ref{tab:yinit} defines two sets (A and B) of three
initial states of the forest and Table~\ref{tab:harv} shows eight alternative 
cases of harvesting schedules.
\begin{table}[!htp]
\caption{
Initial number of trees in size classes. Three alternatives for the initial states of the forest 
denoted by $e$, $n$ and $o$ are considered in two sets A and B.}
\label{tab:yinit}
\begin{center} {\footnotesize
\begin{tabular}{|r|rrr|rrr|}
\multicolumn{7}{c}{}\\ \hline
\multicolumn{1}{|c|}{size}&
\multicolumn{3}{c|}{A}&
\multicolumn{3}{c|}{B}\\
\cline{2-7}
\multicolumn{1}{|c|}{class}&
\multicolumn{1}{c}{e}&
\multicolumn{1}{c}{n}&
\multicolumn{1}{c|}{o}&
\multicolumn{1}{c}{e}&
\multicolumn{1}{c}{n}&
\multicolumn{1}{c|}{o}\\
\hline
 1&  20& 100&   0& 196&1750&  50\\
 2&  20& 100&   0& 162&   0&  25\\
 3&  20& 100&   0& 140&   0&  10\\
 4&  20& 100&   0& 124&   0&   0\\
 5&  20& 100&   0&  75&   0&  25\\
 6&  20&   0&   0&  18&   0& 250\\
 7&  20&   0&   0&   0&   0&  25\\
 8&  20&   0&  50&   0&   0&   0\\
 9&  20&   0&  50&   0&   0&   0\\
10&  20&   0&  50&   0&   0&   0\\
11&  20&   0&  50&   0&   0&   0\\
12&  20&   0&  50&   0&   0&   0\\
\hline
\end{tabular}}
\end{center}
\end{table}
\begin{table}[!htp]
\caption{
Eight cases of harvesting schedules. Interpretation of the 
three character case code xyz is as follows: 
x is l (s) for long (short) transition period,
y is l (s) for long (short) steady state cycle, and
z is d (s) for dense (sparse) harvesting intervals. The length of
transition period is $t^0$ (steps) and $s$ is the steady state cycle (steps).
For time stages in five year intervals, 1 indicates harvesting and 0 no harvesting.}
\label{tab:harv}
\begin{center} {\footnotesize
\begin{tabular}{|c|c|c|llll|}
\multicolumn{7}{c}{}\\
\hline
\multicolumn{1}{|c|}{case}&
\multicolumn{1}{c|}{$t^0$}&
\multicolumn{1}{c|}{$s$}&
\multicolumn{4}{c|}{harvesting stages}\\
\hline
lld& 30& 6&$0~1~0~0~1~0~0~1~0~0$&$1~0~0~1~0~0~1~0~0~1$&$0~0~1~0~0~1~0~0~1~0$&$0~1~0~0~1~0$\\
lls& 30& 6&$0~1~0~0~0~0~0~1~0~0$&$0~0~0~1~0~0~0~0~0~1$&$0~0~0~0~0~1~0~0~0~0$&$0~1~0~0~0~0$\\
lsd& 30& 3&$0~1~0~0~1~0~0~1~0~0$&$1~0~0~1~0~0~1~0~0~1$&$0~0~1~0~0~1~0~0~1~0$&$0~1~0$\\
lss& 30& 3&$0~1~0~0~0~0~0~1~0~0$&$0~0~0~1~0~0~0~0~0~1$&$0~0~0~0~0~1~0~0~0~0$&$0~1~0$\\
sld& 10& 6&$1~0~0~0~1~0~0~1~0~0$&$0~1~0~0~1~0$&&\\
sls& 10& 6&$0~1~0~0~0~0~0~1~0~0$&$0~1~0~0~0~0$&&\\
ssd& 10& 3&$1~0~0~1~0~0~1~0~0~1$&$0~0~1$&&\\
sss& 10& 3&$0~1~0~0~0~0~0~0~1~0$&$0~0~1$&&\\
\hline
\end{tabular}}
\end{center}
\end{table}


Each of the seven methods was run with 1000 independent repetitions 
for each of the $2\times 3\times 8=48$ test problems.
Table~\ref{tab:optval} shows the best objective function value found which we now call the optimum. 
It is interesting to note that for set B of 24 problems 
each of the seven methods found the optimal solution. The same is true
for set A with one exception: for the case e-ssd the optimal value is 14.684
but three of the methods ended up with best value of 14.680. 
Because the difference is only 0.03 \% and given data precision, we may regard 
the value 14.680 optimal as well.

\begin{table}[!htp]
\caption{
Best objective function values found in 7000 runs.}
\label{tab:optval}
\begin{center} {\footnotesize
\begin{tabular}{|c|rrr|rrr|}
\multicolumn{7}{c}{}\\ \hline
\multicolumn{1}{|c|}{case}&
\multicolumn{3}{c|}{A}&
\multicolumn{3}{c|}{B}\\
\cline{2-7}
\multicolumn{1}{|c|}{}&
\multicolumn{1}{c}{e}&
\multicolumn{1}{c}{n}&
\multicolumn{1}{c|}{o}&
\multicolumn{1}{c}{e}&
\multicolumn{1}{c}{n}&
\multicolumn{1}{c|}{o}\\
\hline
lld&  12.544&  8.882&  18.626&  10.475&   8.691&  14.087\\
lls&  12.675&  8.825&  18.784&  10.352&   7.602&  14.221\\
lsd&  12.544&  8.882&  18.626&  10.475&   8.691&  14.087\\
lss&  12.665&  8.815&  18.774&  10.343&   7.593&  14.211\\
sld&  14.749&  8.838&  23.824&  10.488&   7.820&  15.022\\
sls&  12.499&  8.715&  18.581&  10.298&   6.886&  14.065\\
ssd&  14.684&  8.862&  23.767&  10.564&   8.185&  14.961\\
sss&  12.653&  8.745&  18.774&  10.256&   6.705&  14.209\\
\hline
\end{tabular}}
\end{center}
\end{table}
In 34 out of the 48 test problems all seven methods never produced a suboptimal solution
in 1000 trials. Here suboptimal means a local optimum for which the objective function 
value is less than the best found by the method in 1000 trials. 
Among 9 problems, frequency of suboptimal solutions was less than 0.5 \% in 1000 trials.
The general observation is that the likelihood
for suboptimal solution is small. There are a few exceptions concerning set A. 
First, in the case e-ssd all methods most frequently find the second best 
value 14.680 (which lags 0.03 \% from the optimum). Second, in cases e-lld and e-lsd
two of the seven versions produce the second best value 12.537 in more than 10 \% of runs; 
however, the optimum 12.544 is only 0.05 \% higher. 
For the random initialization step in Section~\ref{sec:nonConvex} no suboptimal solution was 
found in 45 out of 48 test cases. In two cases suboptimal frequency was less than 0.5 \% 
and in the case e-ssd only the second best value 14.680 (which lags 0.03 \% from the optimum)
is found in 1000 trials.

Finally, we note that for any particular problem initialization by any of 
the seven methods does not necessarily lead to a local optimum. Due to non-convexities, 
the solver code may stop before optimality conditions are met. In such a case 
we try again and repeat until a successful completion. We observed for the seven 
versions that in 95 \% of the case problems the average number of trials needed 
for success is close to one (the first trial being a success in most cases).
For the initialization in Section~\ref{sec:nonConvex} the average number 
of trials was more than 1.4 in three cases, the maximum average being 2.9 for n-sls.

\section{MIP approach for fitness evaluation}\label{sec:ec2}

In this section we consider mixed integer programing ($MIP$) for 
the problem for which the evolutionary algorithm is proposed in
Sections~\ref{sec:level1}-\ref{sec:nonConvex}.
Let $t^0$ and $t^1$ denote the beginning of the first and second cycle, respectively.
A complicating feature for $MIP$ is that both $t^0$ and $t^1$ are endogenous in the model.
Consequently, for the steady state requirement $x_{t^0}=x_{t^1}$ subscripts of decision 
vectors are endogenous. Furthermore, cash flows $F_t=R_t-C_t$ (refer Equation~\ref{obj}) are treated differently depending 
on whether $t<t^0$, $t^0\leq t<t^1$ or  $t\geq t^1$. We begin this section 
by reformulating the problem as a $MIP$ problem and discuss the implementation of branch-and-bound algorithm
thereafter. Computational results are presented in Section 4.4.1, 
including comparison with the evolutionary algorithm.

Let $T=t_{\max}+s_{\max}$ be an exogenous upper limit for $t^1$.
For all $t=0,1,\dots,T$ let $\delta_t$, $u_t$ and $r_t$ be binary variables such that
$\delta_t=1$ if harvesting takes place at stage $t$ and $\delta_t=0$ otherwise,
$u_t=1$ if $t=t^0$ and $u_t=0$ otherwise, and $r_t=1$ if $t=t^1$ and $r_t=0$ otherwise. 
Hence,
\be
\sum_{t=0}^T u_t =1~~~ {\rm and}~~~ \sum_{t=0}^T r_t =1,
\label{sumuv}
\ee
\be
t^0=\sum_t tu_t ~~~{\rm with} ~~~t_{\min}\leq t^0 \leq t_{\max},
\label{t0}
\ee
\be
t^1=\sum_t tr_t,
\label{t1}
\ee
and the cycle length is
\be
s = t^1-t^0 ~~~{\rm with} ~~~s_{\min}\leq s \leq s_{\max},
\label{s}
\ee
where $s_{\min}$ and $s_{\max}$ are exogenous bounds for cycle length,
and $t_{\min}$ and $t_{\max}$ are exogenous bounds for transient period.

The level of harvest $h_t$ (number of trees at stage $t$ by size class) 
is restricted by
\be
0\leq h_t \leq \delta_t\bar{x}
\label{H}
\ee
where $\bar{x}$ is an exogenous upper bound for the forest state vectors $x_t$ (number of trees 
at stage $t$ by size class). 
Harvests are only allowed during the period $[0,t^1)$. Hence we require
\be
\delta_t \leq \sum_{\tau>t} r_{\tau}.
\label{endH}
\ee
At the beginning of each cycle we require the state to be 
an endogenous vector $x^c$. Hence for all $t\geq t_{\min}$, we require
\be
x_t - \bar{x}(1-u_t-r_t) \leq x^c \leq x_t + \bar{x}(1-u_t-r_t).
\label{yc}
\ee

For all $t$, the cash flow $F_t$ is split into
\be
F_t = F_t^0 + F_t^1
\label{split}
\ee
so that $F_t = F_t^0$ and $F_t^1=0$ for the transient period $t<t^0$, 
$F_t = F_t^1$ and $F_t^0=0$ for the first cycle $t^0\leq t <t^1$, and 
$F_t^0=F_t^1=0$ for $t\geq t^1$. 
This is achieved at an optimum if we require for all $t$ 
\be
-\bar{F}\sum_{\tau> t}u_{\tau} \leq F_t^0 \leq \bar{F}\sum_{\tau> t}u_{\tau} 
\label{c0}
\ee
and
\be
-\bar{F}\sum_{\tau\leq t}(u_{\tau}-r_{\tau})\leq F_t^1 \leq \bar{F}\sum_{\tau\leq t}(u_{\tau}-r_{\tau})
\label{c1}
\ee
where $\bar{F}$ is an exogenous upper bound for cash flows. 

In this notation, the mixed integer programing problem is to find 
length $t^0$ of transition period, 
cycle length $s$, 
forest state $x_t\geq 0$, 
harvesting levels $h_t\geq 0$, 
steady state vector $x^c$, 
cash flows $F_t$, $F_t^0$, $F_t^1$
auxiliary variables, and
binary variables $\delta_t$,  $u_t$,  $r_t$ , for all $t$, to
\be
\max\sum_{t=0}^T \beta^{t\Delta}[F_t^0+F_t^1/(1-\beta^{s\Delta})],
\label{obj_mip}
\ee
subject to (\ref{sumuv})--(\ref{c1}), (\ref{bt})--(\ref{Ct}) and (\ref{xt1})--(\ref{x0}) for all $t$.

We are not aware of any software for solving this problem. Therefore, using AMPL \citep{fourer1990modeling} and 
MINOS \citep{murtagh1978large} as the solver, we implement the well-known branch-and-bound procedure
with two special features. First, to deal with the large number of binary variables, 
$u_t$ and $r_t$ in (\ref{sumuv})--(\ref{s}), we use partitioning. For each node of the tree, there is
a range $[\check{t},\hat{t}]$ for the end $t^0$ of the transition period, and 
a range $[\check{s},\hat{s}]$ for the cycle length $s$. In the root node, $[\check{t},\hat{t}]=[t_{\min},t_{\max}]$
and $[\check{s},\hat{s}]=[s_{\min},s_{\max}]$. Binary variables $u_t$ are fixed to zero for all $t<\check{t}$ and for all
$t>\hat{t}$. For $\check{t}\leq t \leq \hat{t}$, $u_t$ is a continuous variable in the interval $[0,1]$.
Similarly, we fix binary variables $r_t$ to zero if $t<\hat{t}+\check{s}$ or if $t>\hat{t}+\hat{s}$; otherwise
$r_t$ is a continuous variable. Second, after the choice of the node for branching is done, 
if at the optimum of the relaxed problem some of the variables $u_j$ is in the interior of $[0,1]$, 
then branching is based on partitioning the interval $[\check{t},\hat{t}]$ into 
two subsets of equal (or almost equal) cardinality. If no such  variable $u_j$ exists,
then we try similarly to partition $[\check{s},\hat{s}]$. Finally, if all variables $u_t$ and $r_t$
are at levels 0 or 1, then branching is based on one of the binary variables $\delta_t$ in the 
interior of $[0,1]$; if none exists, then stop.

\end{appendices}

\end{document}